\magnification=\magstep1
\font\tensmc=cmcsc10
\hfuzz.3truein
\def\smc{\tensmc}
\def\srm{\sevenrm }

\def\bbox{\quad\hbox{\vrule
 \vbox{\hrule \vskip2pt \hbox{\hskip2pt
 \vbox{\hsize=1pt}\hskip2pt} \vskip2pt\hrule}\vrule}}
\def\qed{$\bbox$}
\def\lessim{\ \lower4pt\hbox{$
 \buildrel{\displaystyle <}\over\sim$}\ }
\def\gessim{\ \lower4pt\hbox{$
 \buildrel{\displaystyle >}\over\sim$}\ }
\def\insim{\ \lower4pt\hbox{$
 \buildrel{\displaystyle \in}\over\sim$}\ }
\def\lessim{\ \lower4pt\hbox{$
 \buildrel{\displaystyle <}\over\sim$}\ }
\def\n{\noindent}

\def\ve{\varepsilon}
\def\vet{\tilde{\varepsilon}}
\def\var{{\rm Var}}
\def\card{{\rm Card}}

\vsize=24.2true cm
     \hsize=15.3true cm
     \nopagenumbers
     \topskip=1truecm
     \headline={\tenrm\hfil\folio\hfil}
     \raggedbottom
     \abovedisplayskip=3 mm % Reduction of space between
                            % text and formulae
     \belowdisplayskip=3 mm %
     \abovedisplayshortskip=0mm
     \belowdisplayshortskip=2 mm
     \normalbaselineskip=12pt %This is default. Do NOT change it!
     \normalbaselines

\hfuzz=.4truein

\topglue.5truein
\centerline{\bf
THE LIL FOR CANONICAL U-STATISTICS OF ORDER 2}\par

\vskip.35truein
\centerline{
{\smc by Evarist Gin\'e$^*$, Stanislaw Kwapie\'n$^{\dagger}$, Rafa{\l}
Lata\l a$^{\dagger}$,
and Joel Zinn}$^\ddagger$}\par
\vskip.08truein\centerline{\it University of Connecticut,
Warsaw University, Warsaw University}
\centerline{\it and Texas A\&M University}\par
\vskip0.7truein
\vfootnote {$\phantom{*}$}{\srm $\scriptstyle{AMS\ 1991\ subject
\ classifications.}$
Primary: 60F15.}
\vfootnote {$\phantom{*}$}{\srm $\scriptstyle{Key\ words\
and\ phrases}$:
$\scriptstyle U$-statistics (canonical or degenerate), law of the iterated
logarithm.}
\vfootnote *{\srm Research partially supported
 by NSF Grant No. DMS-96--25457.}
\vfootnote \dag{\srm Research partially supported by Polish Grant KBN 2
PO3A 043
15.}
\vfootnote \ddag{\srm Research partially supported by NSF Grant No.
DMS-96-26778.}
\centerline{\bf Abstract} \vskip0.15truein
{\narrower{\srm Let $\scriptstyle {X,\ X_i,\ i\in\bf
N}$, be independent identically distributed random variables and let
$\scriptstyle h(x,y)=h(y,x)$ be a measurable function of two variables.
It is shown that the bounded law of the iterated logarithm, $\scriptstyle
\limsup_n (n\log\log n)^{-1}\big|\sum_{1\le i< j\le n}h(X_i,X_j)\big|<\infty$
a.s., holds if and only if the following three conditions are satisfied:
$\scriptstyle h$ is canonical for the law of
$\scriptstyle X$ (that is,
$\scriptstyle{E}h(X,y)=0$ for almost all $\scriptstyle y$) and there
exists
$\scriptstyle C<\infty$ such that, both,
$\scriptstyle E(h^2(X_1,X_2)\wedge u)\le C\log\log u$ for all large
$\scriptstyle u$ and
$\scriptstyle \sup\{
Eh(X_1,X_2)f(X_1)g(X_2):\|f(X)\|_2\le1,\|g(X)\|_2\le1, \|f\|_\infty<\infty,
\|g\|_\infty<\infty\}\le C$.
\par}}  \vskip1.1truein
\centerline{May 1999}
\vfill\eject

\n{\bf 1. Introduction.} Although $U$--statistics
(Halmos, 1946; Hoeffding, 1948) are relatively simple probabilistic objects,
namely averages over an i.i.d. sample
$X_1,\dots, X_n$ of measurable functions (kernels) $h(x_1,\dots,x_m)$ of
several
variables, their asymptotic theory is only recently attaining a satisfactory
degree of completeness: see e.g. Rubin and Vitale (1980), Gin\'e and Zinn
(1994), Zhang (1999) and Lata\l a and Zinn (1999) on necessary and sufficient
conditions for the central limit theorem and the law of large numbers. We are
interested here in the law of the iterated logarithm for $U$-statistics
based on
canonical (or completely degenerate) kernels, that is, on kernels whose
conditional expectation given any $m-1$ variables is zero, and only for $m=2$.

$U$-statistics with nondegenerate
kernels behave, as is well known, like sums of independent random
variables, and
the LIL in this case was proved by Serfling (1971). The LIL for canonical (or
completely degenerate) kernels $h$ with finite absolute moment of order
$2+\delta$,
$\delta>0$, was obtained by Dehling, Denker and Philipp (1984, 1986), and with
finite second moment by Dehling (1989) and Arcones and Gin\'e (1995). Gin\'e
and Zhang (1996) showed that there exist degenerate kernels $h$ with infinite
second moment such that, nevertheless, the corresponding
$U$-statistics satisfy the law of the iterated logarithm, and obtained a
necessary
integrability condition as well. This last article and Goodman's (1996) also
contain LIL's under assumptions that do not imply finiteness of the second
moment of
$h$, but that fall quite short from being necessary. The LIL
for finite sums of products $\sum_{i=1}^k\lambda_i\phi_i(x_1)\cdots\phi_i(x_m)$
is easier ($E h^2<\infty$ is necessary) and was considered by Teicher
(1995) for $k=1$ and by  Gin\'e and Zhang (1996) for any $k<\infty$.
In the present
article the bounded LIL problem is solved for kernels of order 2. Next we
describe our result and comment on its (relatively involved) proof.

In what follows, $X, X_i$, $i\in\bf N$, are independent identically
distributed random variables taking values on some measurable space $(S,{\cal
S})$, and
$h:S^2\mapsto\bf R$ is a measurable function that we assume, without loss of
generality (for our purposes), symmetric in its entries, that is,
$h(x,y)=h(y,x)$ for all $x,y\in S$. When $h$ is integrable we say that it is
 canonical, or degenerate, for the law of $X$ if $E h(X,y)=0$ for almost all
$y\in S$ (relative to the law of $X$). The natural LIL normalization for
$U$-statistics corresponding to degenerate kernels of order 2 is $n\log\log n$
as is seen with the following example. A simple canonical kernel for
$S=\bf R$ and $X$ integrable with
$E X=0$ is
$h(x,y)=xy$. For this example, if moreover $E X^2<\infty$ then, by the LIL
and the law of large numbers for sums of independent random variables, we have
$$\limsup_n{1\over 2n\log\log n}\bigg|\sum_{i\ne j\le
n}X_iX_j\bigg|=\limsup_n\biggl[{1\over\sqrt{2n\log\log
n}}\sum_{i=1}^nX_i\biggr]^2 =\var X.$$

Our main result is as follows:

\proclaim{\smc Theorem 1.1}. Let $X,Y,X_i$, $i\in\bf N$, be i.i.d. random
variables taking values in
$(S,{\cal S})$ and let
$h:S^2\mapsto \bf R$ be a measurable function of two variables. Then,
$$\limsup_n{1\over n\log\log n}\bigg|\sum_{1\le i\ne j\le
n}h(X_i,X_j)\bigg|<\infty\ \ {\rm a.s.}\eqno(1.1)$$
if an only if the following three conditions hold: \hfil\break
a) $h$ is canonical for the law of $X$\hfill\break
and there exists $C<\infty$ such that \hfil\break
b)  for all $u\ge10$,
$$E(h^2(X,Y)\wedge u)\le C\log\log u,\eqno(1.2)$$
and\hfil\break
c) $$\eqalignno{\sup\bigl\{E h(X,Y)f(X)g(Y):E f^2&(X)\le
1,E g^2(X)\le 1,&\cr
&
\|f\|_\infty<\infty,\|g\|_\infty<\infty\bigr\}\le C.&(1.3)\cr}$$

It is easily seen that condition b) implies
$$E {h^2\over(\log\log(|h|\vee e^e)^{1+\delta}}<\infty\eqno(1.4)$$
for all $\delta >0$ (and is implied by $E
h^2/\log\log(|h|\vee e^e)<\infty$. In particular condition b) ensures the
existence of the integrals in conditions a) and c). Condition c) implies that
the operator defined on
$L_\infty({\cal L}(X))$ by $Hf(y)=E h(X,y)f(X)$ takes values in $L_2({\cal
L}(X))$ and extends as a bounded operator to all of
$L_2({\cal L}(X))$. Moreover, if with a slight abuse of notation we set
$E_X h(X,Y)f(X):=Hf(Y)$ for $f\in L_2$, then  condition b) is
equivalent to
$$E_Y\bigl(E_X h(X,Y)f(X)\bigr)^2\le C^2E f^2(X)\ \ {\rm for\
all}\ f\in L_2.\eqno(1.5)$$
(Here and in what follows, $E_X$ (resp. $E_Y$) indicates
expectation with respect to $X$ (resp. $Y$) only.)

The integrability condition b) was proved to be necessary for the LIL (1.1) by
Gin\'e and Zhang (1996), whereas the idea for condition c) comes from
Dehling (1989) who showed that if
$h(x,y)$ is canonical and square integrable then
$$\eqalign{{\rm lim\ set}&\biggl\{{1\over2n\log\log n}\sum_{1\le i\ne j\le
n}h(X_i,X_j)\biggr\}\cr
&=
\bigl\{E h(X,Y)f(X)f(Y):E f^2(X)\le 1\bigr\}\ \ {\rm a.s.}\cr}$$

We will not prove Theorem 1.1 directly, but instead we will prove first that
conditions b) and c) are necessary and sufficient for a decoupled and
randomized
version of the LIL, namely, for
$$\limsup_n{1\over n\log\log n}\bigg|\sum_{1\le i, j\le
n}\ve_i\vet_jh(X_i,Y_j)\bigg|<\infty\ \ {\rm a.s.},\eqno(1.6)$$
where $\{\ve_i\}$ is a Rademacher sequence
independent of all the other variables. (We recall that a Rademacher
sequence is
a sequence of independent random variables taking on only the values $1$ and
$-1$, each with probability 1/2.) The reasons for this are multiple. One
is that necessity of condition c) follows as a consequence of a recent
result of
Lata\l a (1999) on estimation of tail probabilities of Rademacher chaos
variables. Another reason is that, because of the Rademacher multipliers,
truncation of the kernel will result in symmetric, and hence mean
zero, variables; this is important since the proof of sufficiency contains
several relatively complicated truncations of $h$.  Moreover, part of the core
of the proof of sufficiency consists of an iterative application of an
exponential bound for sums of independent random variables and vectors, and
having decoupled expressions makes this iteration possible (although we could
use, alternatively, an exponential inequality for martingale differences that
does not require decoupled expressions).

The exponential inequality in question is Talagrand's (1996) uniform Prohorov
inequality. This inequality depends on two parameters, the $L_\infty$ bound of
the variables and the weak variance of their sum, and to apply it iteratively
requires not only that
$h$ be truncated at a low level, but that the conditional second moments of
these
truncations of $h$ be small as well. This explains the relatively complicated
multi-step truncation procedure in the proof of sufficiency.

Finally, the limit (1.6) will
imply the limit (1.1) by a two stage symmetrization argument that will also
require control of the conditional expectations of the sums; this control
will be
achieved once more, again after multiple truncations, by means of Talagrand's
exponential inequality.

Section 2 contains several known results needed in the sequel. Section 3 is
devoted to the proof of the LIL for decoupled, randomized kernels, and Section
4 reduces the LIL for canonical kernels to this case. In Section 5 we complete
the proof of Theorem 1.1 and make several comments about the limsup in (1.1)
and the limit set of the LIL sequence.

We adhere in what follows to the following notation (some of it already set up
above):
\item{$\diamond$} $h$ is a measurable real function of two variables defined on
$(S^2,{\cal S}\otimes{\cal S})$, symmetric in its entries.

\item{$\diamond$} $X,X_{1},X_{2},\ldots$ and $Y,Y_{1},Y_{2},\ldots$ denote two
independent, equidistributed sequences of i.i.d. $S$-valued random variables.

\item{$\diamond$} We write $Ef(h)$ for $Ef(h(X,Y))$, and $E_X$, $\Pr_X$
(resp.
$E_Y$, $\Pr_Y$) denote expected value and probability with respect to the
random variables $X,X_i$ (resp. $Y,Y_i$) only.

\item{$\diamond$} $\ve_{1},\ve_{2},\dots,$ and $\vet_{1},\vet_{2},\ldots$  are
two independent Rademacher sequences,
 independent of all other random variables.

\item{$\diamond$}  We write $L_{2}x$ and $L_{3}x$
 instead of $L(L(x))$ and $L(L(L(x)))$, where $L(x)=\max(\log x,1)$.

\item{$\diamond$} In all proofs $\tilde{C}$ denotes a universal constant
which
  may change from line to line but does not depend on any parameters.

\vskip.15truein
\noindent {\bf 2. Preliminary results.} For convenience, we isolate in this
section several known results needed below.
\vskip.1truein

\n {\it (A) Hoeffding's decomposition.} The $U$-statistics with kernel $h$ (not
necessarily symmetric in its entries) based on $\{X_i\}$ are defined as
$$U_n(h)={1\over n(n-1)}\sum_{1\le i\ne j\le n}h(X_i,X_j),~~~n\in{\bf N}.$$
By considering instead the kernel
$\tilde{h}(x,y)=\bigl(h(x,y)+h(y,x)\bigr)/2$, we have
$$U_n(h)=U_n(\tilde{h})={1\over n(n-1)}\sum_{1\le i\ne j\le
n}\tilde{h}(X_i,X_j)={n\choose 2}^{-1}\sum_{1\le i< j\le n}\tilde h(X_i,X_j).$$
So, {\it we will assume $h$ symmetric in its entries} in all that follows.

Suppose
$E|h(X,Y)|<\infty$. Then,
$$\eqalignno{h(x,y)-Eh(X,Y)&=\bigl[h(x,y)-E_Yh(x,Y)-E_Xh(X,y)+Eh(X,Y)\bigr]&\cr
&~~~~~~~+
\bigl[E_Yh(x,Y)-Eh(X,Y)\bigr]+\bigl[E_Xh(X,y)-Eh(X,Y)\bigr]&\cr
&:=\pi_2h(x,y)+\pi_1h(x)+\pi_1h(y),&(2.1)\cr}$$
where the identities hold a.s. for ${\cal L}(X)\times{\cal L}(X)$. The kernel
$\pi_2h$ is canonical (or degenerate) for the law of $X$ as $E_X\pi_2h(X,Y)=
E_Y\pi_2h(X,Y)=0$ a.s., and $\pi_1h(X)$ is centered. This decomposition of
$h$ gives rise to {\it Hoeffding's decomposition} of the corresponding
$U$-statistics,
$$\sum_{1\le i<j\le n}h(X_i,X_j)=\sum_{1\le i< j\le
n}\pi_2h(X_i,X_j)+(n-1)\sum_{i=1}^n\pi_1h(X_i)+{n\choose 2}Eh(X,Y),\eqno(2.2)$$
and of their decoupled versions,
$$\eqalignno{\sum_{1\le i, j\le n}h(X_i,Y_j)=\sum_{1\le i, j\le
n}&\pi_2h(X_i,Y_j)
+n\sum_{i=1}^n\pi_1h(X_i)&\cr
&+n\sum_{i=1}^n\pi_1h(Y_i)
+n^2Eh(X,Y).&(2.3)\cr}$$

\vskip.1truein
\n{\it (B) The equivalence of several LIL statements.} The following lemma
contains necessary randomization and integrability conditions for the LIL:

\proclaim{\smc Lemma 2.1}. {\rm (Gin\'e and Zhang, 1996)}. (a) ({\rm
Integrability.}) There exists a universal constant $K$ such that, if
$$\sum_{n=1}^\infty\Pr\biggl\{{1\over 2^nLn}\Big|\sum_{1\le i,j\le
2^n}\ve_i\vet_jh(X_i,Y_j)\Big|>C\biggr\}<\infty\eqno(2.4)$$
for some $C<\infty$, then
$$\limsup_{u\to\infty}{E\bigl(h^2(X,Y)\wedge u\bigr)\over L_2u}\le
KC^2.\eqno(2.5)$$
(b) ({\rm Randomization and decoupling, partial.}) The LIL
$$\limsup_n{1\over nL_2n}\bigg|\sum_{1\le i<j\le n}h(X_i,X_j)\bigg|\le C\ \
{\rm a.s.}\eqno(2.6)$$
for some $C<\infty$ implies
$$\sum_{n=1}^\infty\Pr\biggl\{{1\over 2^nLn}\max_{k\le 2^n}\Big|\sum_{1\le
i,j\le
k}\ve_i\vet_jh(X_i,Y_j)\Big|>2^7C\biggr\}<\infty.$$
In particular, the LIL implies both the
integrability condition (2.5) and the randomized and decoupled LIL, that is,
$$\limsup_n{1\over nL_2 n}\bigg|\sum_{1\le i, j\le
n}\ve_i\vet_jh(X_i,Y_j)\bigg|\le D\ \ {\rm a.s.}\eqno(2.7)$$
with $D=KC$ for some universal constant $K$.

Part (a) is contained in the proof of Theorem 3.1 in Gin\'e and Zhang (1996),
while part (b) is the content of Theorem 3.1 and Lemma 3.3 there.

We recall that the limsups at the left hand sides of (2.6) and (2.7) are
always a.s. constant (finite or infinite) by the Hewitt-Savage zero-one law.

Decoupling gives the following equivalence between the LIL and its decoupled
version.

\proclaim{\smc Lemma 2.2}. (a) The LIL (2.6) is equivalent to the
decoupled LIL, that is, to
$$\limsup_n{1\over nL_2n}\Big|\sum_{1\le i\ne j\le n}h(X_i,Y_j)\Big|\le D\ \
{\rm a.s.}\eqno(2.8)$$
for some $D<\infty$, meaning that if (2.6) holds for $C$ then (2.8) holds for
$D=KC$ and that if (2.8) holds for $D$ then (2.6) holds for $C=KD$, where $K$
is a universal constant.\hfil\break
(b) The decoupled and randomized LIL (2.7) is
equivalent to the randomized LIL
$$\limsup_n{1\over nL_2 n}\bigg|\sum_{1\le i\ne j\le
n}\ve_i\ve_jh(X_i,X_j)\bigg|\le C\ \ {\rm a.s.}\eqno(2.9)$$
for some $C$ finite (with  $C$ and $D$ related as in part (a)).\hfil\break
(c) The LIL (2.7) implies convergence of the series (2.4) for some
$C=KD<\infty$, $K$ a universal constant, hence it also implies the
integrability
condition (2.5) (with $C$ replaced by $D$).

\n{\smc Proof}. (a) We can equivalently write (2.6) as
$$\lim_{k\to\infty} \Pr\biggl\{\sup_{n\ge k}{1\over nL_2n}\Big|\sum_{1\le i\ne
j\le n}h(X_i,X_j)\Big|\ge C\biggr\}=0$$
for some $C<\infty$, hence as
$$\lim_{k\to\infty} \Pr\biggl\{\Big\|\sum_{1\le i\ne
j<\infty}h_{i\vee j,k}(X_i,X_j)\Big\|\ge C\biggr\}=0,$$
where
$$h_{i,k}:=\biggl({h\over kL_2k},{h\over(k+1)L_2(k+1)},\dots,{h\over
nL_2n},\dots\biggr)$$
if $i\le k$ and
$$h_{i,k}:=\biggl(0,{\buildrel i-k\over\dots},0,{h\over
iL_2i},{h\over(i+1)L_2(i+1)},\dots,{h\over nL_2n},\dots\biggr)$$
if $i>k$ are $\ell_\infty$-valued functions and $\|\cdot\|$ denotes the sup of
the coordinates. Then, the decoupling inequalities of de la Pe\~na and
Montgomery-Smith (1994) apply to show that the above tail probabilities are
equivalent up to constants to those of the corresponding decoupled expressions,
thus giving the equivalence between (2.6) and (2.8).

(b) If (2.9) holds, then (2.7) without  diagonal terms (that is, without the
summands corresponding to
$i=j$) holds too by the first part of the proof applied to the kernel
$\alpha\beta h(x,y)$. Moreover, (2.9) implies the integrability condition (2.5)
by Lemma 2.1 (note that if $\{\varepsilon_i^{(j)}\}$,
$j=1,2,3$, are three independent Rademacher sequences, then
$\{\varepsilon_i^{(1)}\varepsilon_i^{(2)}\}$ and
$\{\varepsilon_i^{(1)}\varepsilon_i^{(3)}\}$ are also independent Rademacher
sequences) and, as a consequence, $h$ is integrable. Hence, by the law of large
numbers, the diagonal in (2.7) is irrelevant, showing that (2.7) holds with the
diagonal included. If (2.7) holds, then we also have
$E|h|<\infty$: a modification of the proof of the converse central limit
theorem in Gin\'e and Zinn (1994), consisting in replacing use of the law of
large numbers by use of inequality (3.7) in Gin\'e and Zhang (1996), shows that
if the sequence
$\bigl\{(nL_2n)^{-1}\sum_{i,j\le n}\ve_i\vet_jh(X_i,Y_j)\}\bigr\}$ is
stochastically bounded, then $Eh^2(X,Y)\wedge u\le C (L_2u)^2$ for some
$C<\infty$, in particular, that $E|h|<\infty$. So, we can delete the
diagonal in
(2.7), and then apply the first part of the lemma to undo the decoupling.

(c) Statement (c) follows from (b) because, by Lemma 2.1, (2.9) implies
convegence of the series (2.4) for some $C<\infty$.
\qed
\vskip.05truein

The following lemma, together with the previous ones, will
allow blocking and will reduce the proof of sufficiency of the LIL to showing
that a series of tail probabilities converges (just as with sums of i.i.d
random variables).

\proclaim{\smc Lemma 2.3}. There exists a universal constant $C<\infty$
such that
for any kernel
$h$
 and any two sequences $X_{i}$, $Y_{j}$ of i.i.d. random variables we have
$$ \Pr\biggl\{\max_{k\leq m,l\leq n}\Big|\sum_{i\leq k,j\leq
l}h(X_{i},Y_{j})\Big|
    \geq t\biggr\}\leq
   C\Pr\bigg\{\Big|\sum_{i\leq m,j\leq n}h(X_{i},Y_{j})\Big|\geq
t/C\bigg\}\eqno(2.10)$$
 for all $m,n\in\bf N$ and for all $t>0$.

\n{\smc Proof.} Montgomery-Smith's (1993) maximal inequality for i.i.d. sums
asserts that if
$Z_{i}$ are i.i.d.\
 r.v.'s with values in some Banach space $B$ then for some universal constant
 $C_{1}$ and all $t>0$ we have
$$\Pr\bigg\{\max_{k\leq m}\Bigl\|\sum_{i\leq k}Z_{i}\bigr\|\geq t\bigg\}\leq
   C_{1}\Pr\bigg\{\Bigl\|\sum_{i\leq m}Z_{i}\Bigr\|\geq t/C_{1}\biggr\}.$$
 We apply this inequality to $B=\ell^{n}_{\infty}$ and
 $Z_{i}=\bigl(\sum_{j\leq l}h(X_{i},y_{j}):l\leq n\bigr)$ for fixed values of
 $y_{1},\ldots,y_{n}$ to get
$$\Pr\biggl\{\max_{k\leq m,l\leq n}\Big|\sum_{i\leq k,j\leq l}h(X_{i},Y_{j})
   \Big|\geq t\biggr\}\leq
   C_{1}\Pr\biggl\{\max_{l\leq n}\Bigl|\sum_{i\leq m,j\leq
l}h(X_{i},Y_{j})\Bigr|
   \geq t/C_{1}\biggr\}.$$
 In a similar way we may prove
$$\Pr\biggl\{\max_{l\leq n}\Bigl|\sum_{i\leq m,j\leq l}h(X_{i},Y_{j})\Bigr|
   \geq t/C_{1}\bigg\}\leq
   C_{1}\Pr\biggl\{\sum_{i\leq m,j\leq n}\Bigl|h(X_{i},Y_{j})\Bigr|\geq
t/C_{1}^{2}
   \biggl\}.$$
  Thus the assertion holds with $C=C_{1}^{2}$. \qed

\vskip.05truein

\proclaim{\smc Corollary 2.4}. If
$$\sum_{n=1}^\infty\Pr\biggl\{{1\over 2^nLn}\Big|\sum_{1\le i,j\le
2^n}h(X_i,Y_j)\Big|>C\biggr\}<\infty\ \ {\rm a.s.}\eqno(2.11)$$
for some $C<\infty$, then there is a universal constant $K$ such that
$$\limsup_n{1\over nL_2n}\Big|\sum_{1\le i,j\le
n}h(X_i,Y_j)\Big|\le KC\ \ {\rm a.s.}\eqno(2.12)$$

\n{\smc Proof}. Since, for any $0<D<\infty$,
$$\eqalign{\Pr\biggl\{\sup_{n\ge N}&{1\over nL_2n}\Big|\sum_{1\le i,j\le
n}h(X_i,Y_j)\Big|>D\biggr\}\cr&\le
\Pr\biggl\{\sup_{k>[\log N/\log 2]}\max_{2^{k-1}\le n\le2^k}{3\over
2^kLk}\Big|\sum_{1\le i,j\le n}h(X_i,Y_j)\Big|>D\biggr\}\cr
&\le\sum_{k>[\log N/\log 2]}\Pr\biggl\{\max_{2^{k-1}\le n\le2^k}\Big|\sum_{1\le
i,j\le n}h(X_i,Y_j)\Big|>{D2^kLk\over 3}\biggr\},\cr}$$
the result follows from Lemma 2.3. \qed

\vskip.05truein

Applying Corollary 2.4 to the kernel $\alpha\beta h(x,y)$ we obtain the
converse of Lemma 2.2(c). Hence,

\proclaim{\smc Corollary 2.5}. Consider the statements
$$\limsup_n{1\over nL_2 n}\bigg|\sum_{1\le i, j\le
n}\ve_i\vet_jh(X_i,Y_j)\bigg|\le C\ \ {\rm a.s.}$$
and
$$\sum_{n=1}^\infty\Pr\biggl\{{1\over 2^nLn}\Big|\sum_{1\le i,j\le
2^n}\ve_i\vet_jh(X_i,Y_j)\Big|>D\biggr\}<\infty.$$
There is a universal constant $K$ such that if the first statement holds for
some $C<\infty$ then the second holds for $D=KC$, and conversely, if the second
holds for some $D<\infty$ then so does the first, for $C=KD$.

We will also require the following partial converse to Lemma 2.1(b) regarding
the regular LIL and convergence of series of tail probabilities:

\proclaim{\smc Corollary 2.6}. Suppose $E|h|<\infty$. If
$$\sum_{n=1}^\infty\Pr\biggl\{{1\over 2^nLn}\Big|\sum_{1\le i,j\le
2^n}h(X_i,Y_j)\Big|>C\biggr\}<\infty\ \ {\rm a.s.}$$
for some $C<\infty$
then the LIL holds, that is, there is a universal constant $K$ such that
$$\limsup_n{1\over nL_2n}\bigg|\sum_{1\le i<j\le n}h(X_i,X_j)\bigg|\le KC\ \
{\rm a.s.}$$

\n{\smc Proof}. Convergence of the series implies (2.12), that is,  the
decoupled
LIL with diagonal terms included. Since $E|h|<\infty$, the diagonal terms are
irrelevant and therefore the decoupled LIL (2.8) holds. The result now follows
from Lemma 2.2(a). \qed

In Section 4 we will apply the conclusion of Corollary 2.6 under the assumption
that the decoupled and randomized LIL (2.7) holds: this is possible because
(2.7) implies integrability of
$h$, as indicated in the proof of Lemma 2.2(b).

\vskip.1truein
\n{\it (C) Inequalities.} As mentioned in the Introduction, the following two
inequalities will play a basic role in the proof of Theorem 1.1. The first
consists of a sharp estimate of the tail probabilities of Rademacher chaos
variables  (it is in fact part of a sharper two sided estimate).

\proclaim{\smc Lemma 2.7}. {\rm (Lata\l a, 1999).} There exists a universal
constant $c>0$ such that, for all matrices $(a_{i,j})$ and for all $t>0$,
$$\Pr\biggl\{\Big|\sum_{i,j}a_{i,j}\ve_i\vet_j\Big|\ge
c|||(a_{i,j})|||_t\biggr\}\ge c\wedge e^{-t},\eqno(2.13)$$
where $|||(a_{i,j})|||_t$ is defined as
$$|||(a_{i,j})|||_t:=\sup\Bigl\{\sum_{i,j}a_{i,j}b_ic_j:\sum_ib_i^2\le t,
\sum_j
c_j^2\le t, |b_i|,|c_j|\le 1\ {\rm for\ all}\ i,j\Bigr\}.\eqno(2.14)$$

The second is a uniform Prohorov inequality due to Talagrand. It combines
Theorem 1.4 in Talagrand (1996) with Corollary 3.4 in Talagrand (1994).

\proclaim{\smc Lemma 2.8}. {\rm (Talagrand, 1996)}. Let $\{X_{i}\}$,
$i=1,\dots,n$ for any $n\in\bf N$, be
 independent random variables with values in a measurable space
 $(S,{\cal S})$, let $\cal F$ be a countable class of measurable
 functions on $S$ and let
$$Z:=\sup_{f\in{\cal F}}\sum_{i=1}^nf(X_i).$$
 There exists a universal constant $K$ such that for
all
$t>0$ and $n\in\bf N$, if
$$\max_{1\le i\le n}\sup_{f\in{\cal F}}{\rm
ess~sup}_{\omega\in\Omega}\big|f(X_i(\omega))\big|\le U,\ \
E\Bigl(\sup_{f\in{\cal F}}\sum_{i=1}^nf^2({\bf X}_i)\Bigr)\le V$$
and $$\sup_{f\in{\cal F}}\sum_{i=1}^nEf^2(X_i)\le\sigma^2,$$
then
 $$\eqalignno{\Pr\Bigl\{|Z-E Z|\geq t\Bigr\}&\leq K\exp\biggl(-{t\over KU}
   \log\Big(1+{tU\over V}\Bigr)\biggr)&\cr
&\le K\exp\biggl(-{t\over KU}
   \log\Big(1+{tU\over \sigma^2+8UE|Z|}\Bigr)\biggr).&(2.15)\cr}$$

In fact, we will only use the corresponding deviation inequality, that is, the
bound (2.5) for $\Pr\{Z>EZ+t\}$. Ledoux (1987) contains a simple proof of this
result based on logarithmic Sobolev inequalities.

When $\cal F$ consists of a single function $f$ and the variables
$f(X_i)$  are centered this inequality reduces, modulo constants, to the
classical Prohorov inequality. For convenience, we will refer below to
Lemma 2.8
even in cases when Prohorov's inequality suffices.

\vskip.1truein
\n{\bf 3. Symmetrized kernels.} In this section we prove the following
theorem,
which constitutes
the basic component of the proof of Theorem 1.1.

\proclaim{\smc Theorem 3.1}. The decoupled and randomized LIL holds,  that is,
$$\limsup_n{1\over n\log\log n}\Bigl|\sum_{1\le i,j\leq n}
   \ve_{i}\vet_{j}h(X_{i},Y_{j})\Bigr|<\infty \ \ {\rm a.s.}\eqno(3.1)$$
if and only if the following two conditions
 are satisfied for some $C<\infty$:
$$E\min(h^{2},u)\leq CL_{2}u \ \ {\rm for\ all}\ u>0,\eqno(3.2)$$
and
$$\eqalignno{\sup\bigl\{Eh(X,Y)f(X)g(Y):Ef^{2}(X)&\leq 1,
Eg^{2}(Y)\leq 1,&\cr
&\|f\|_\infty<\infty,\|g\|_\infty<\infty\bigr\}\leq C<\infty.
&(3.3)\cr}$$

\n{\smc Remark.}  We recall that, by Corollary 2.5, a necessary and
sufficient condition for the LIL (3.1)  to hold is that
$$\sum_{n=1}^\infty\Pr\biggl\{{1\over 2^nLn}\Big|\sum_{1\le i,j\le
2^n}\ve_i\vet_jh(X_i,Y_j)\Big|>C\biggr\}<\infty\eqno(3.4)$$
for some $C<\infty$.

\vskip.1truein
\n{\smc Proof of necessity.} The integrability condition (3.2) is necessary
for (3.1) by Lemma 2.2(c). The necessity of (3.3) will follow from Lemma 2.7.
For this, we estimate first $\big|\big|\big|\bigl(h(X_i,Y_j):i,j\le
2^n\bigr)\big|\big|\big|_{\log n}$, where $|||\cdot|||_t$ is as defined in
(2.13).  Suppose that
$f,g\in L_{\infty}$ are such that
$Ef^{2}(X)=Eg^{2}(X)=1$ and set
$$K:=|Eh(X,Y)f(X)g(Y)|,\eqno(3.5)$$
that we can assume strictly positive. Note that the integral exists by
(3.2). Then by the SLLN for i.i.d. r.v.'s and
$U$-statistics we have a.s.
$$n^{-1}\sum_{i\leq n}f^{2}(X_{i})\rightarrow Ef^{2}=1,\ \
n^{-1}\sum_{j\leq n}g^{2}(Y_{j})\rightarrow Eg^{2}=1$$
and
$$n^{-2}\Big|\sum_{i,j\leq n}h(X_{i},Y_{j})f(X_{i})g(Y_{j})\Big|\rightarrow
 \big|Eh(X,Y)f(X)g(Y)\big|.$$
So, for large enough $n$,
$$\Pr\biggl\{2^{-n}\sum_{i\leq 2^{n}}f^{2}(X_{i})\leq 2\biggr\}\geq
{3\over4},\ \ \
\Pr\biggl\{2^{-n}\sum_{j\leq 2^{n}}g^{2}(Y_{j})\leq 2\biggr\}\geq
{3\over4}$$     and
$$\Pr\biggl\{2^{-2n}\Bigl|\sum_{i,j\leq 2^{n}}h(X_{i},Y_{j})f(X_{i})g(Y_{j})
  \Bigr|\geq K/2\biggr\}\geq {3\over4}$$
with $K$ as in (3.5). Since $f,g\in L_{\infty}$ we have that, for large enough
$n$,
$$\biggl|\sqrt{{\log n\over2^{n+1}}}f(X_{i})\biggr|,
  \biggl|\sqrt{{\log n\over2^{n+1}}}g(Y_{j})\biggr|\leq 1 \ \ {\rm a.s. }$$
Then, it follows directly from the definition of $|||\cdot|||_t$ that, on the
intersection of the above five events, we have the bound
$$\big|\big|\big|\bigl(h(X_{i},Y_{j}):i,j\leq
2^{n}\bigr)\big|\big|\big|_{\log n}\ge
   K2^{n-2}\log n.$$
Therefore, for large $n$,
$$\Pr\Bigl\{\big|\big|\big|\bigl(h(X_{i},Y_{j}):i,j\leq
2^{n}\bigr)\big|\big|\big|_{\log n}\geq
   K2^{n-2}\log n \Bigr\}\geq {1\over4}.$$
Then, Lemma 2.7 implies that, for all $n$ large enough,
$$\Pr\biggr\{\Bigr|\sum_{i,j\leq
2^{n}}h(X_{i},Y_{j})\ve_{i}\vet_{j}\Bigl|\geq
   cK2^{n-2}\log n\biggl\}\geq {1\over 4}e^{-\log n}= {1\over 4n}.$$
By (3.4), this implies that if the LIL holds then $K$ is uniformly bounded,
proving necessity of condition (3.3). \qed

\vskip.05truein

Before starting the proof of sufficiency, it is convenient to show how the
integrability condition (3.2) limits the sizes of certain truncated conditional
second moments. To simplify notation, we define
$$f_{n}(x)=E_{Y}\min\bigl(h^{2}(x,Y),2^{4n}\bigr)\ \ {\rm
and}\ \ f_{n}(y)=E_{X}\min\bigl(h^{2}(X,y),2^{4n}\bigr).\eqno(3.6)$$

\proclaim{\smc Lemma 3.2}. For any kernel $h$ satisfying condition (3.2)
we have that, for all
  $a>0$,
$$\sum_{n}2^{n}{\Pr}_{X}\Bigl\{E_{Y}\min\bigl(h^{2}(X,Y),2^{an}\bigr)\geq
    2^{n}(\log n)^2\Bigr\}<\infty.\eqno(3.7)$$
Moreover,
  $$ \sum_{n}{2^{n}\over(\log n)^{k}}\Pr\Bigl\{f_{n}(X)\geq 2^{n}(\log n)^{2-k}
    \Bigr\}<\infty ~~for\ all\ k\geq 0.\eqno(3.8)$$

\n{\smc Proof.} For $a$ fixed, we set $\gamma_{k}=\exp(2^{k+1})$ and
$\tilde{f}_{k}(X)=E_{Y}\min(h^{2},2^{a\gamma_{k}})$. Then,
$$\eqalignno{\sum_{2^{k}\leq \log n\leq 2^{k+1}}2^{n}{\Pr}_{X}\Bigl\{
  E_{Y}\min\bigl(h^{2}(X,Y)&,2^{an}\bigr)\geq 2^{n}(\log n)^2\Bigr\}&\cr
&\leq
\sum_{2^{k}\leq \log n\leq 2^{k+1}}
  2^{n}{\Pr}_{X}\bigl\{\tilde{f}_{k}(X)\geq 2^{n+2k}\bigr\}&\cr
&\leq
  E\sum_{n}2^{n}I\bigl(\tilde{f}_{k}(X)\geq 2^{n+2k}\bigr)&\cr
&\leq
2^{1-2k}E\tilde{f}_{k}(X)\leq
2^{1-2k}CL_{2}(2^{a\gamma_{k}})&\cr
&\leq 2^{1-2k}C(\log a+2^{k+1}).&(3.9)\cr}$$
Convergence in (3.7) follows from (3.9). Condition (3.8) is an
easy consequence of (3.7) (as can be seen e.g. by making the
approximate change of variables $2^n/(\log n)^k\simeq 2^m$ in (3.8) and
comparing
with (3.7) for $a>4$).
\qed

\vskip.1truein
\n{\smc Proof of sufficiency.} Since this is only a matter of normalization
we will assume
that conditions (3.2) and (3.3) are satisfied with $C=1$.
By the Remark below Theorem 3.1, proving the LIL is equivalent to showing that
the series (3.4) converges for some
$C<\infty$.
To establish this we will show in several steps that we may suitably truncate
$h$ by proving inequalities of the form
$$
\sum_{n}\Pr\biggl\{\Bigl|\sum_{i,j\le
2^{n}}\ve_{i}\vet_{j}h_{n}(X_{i},Y_{j})\Bigr|
  \geq C2^{n}\log n\biggr\}<\infty,\eqno(3.10)
$$
where $h_{n}:=hI_{A_{n}}$ and $A_{n}$ are suitably chosen subsets of
the product space. Then, we
will apply Lemma 2.8 conditionally to the truncated $h$ (several times, and
after some additional preparation).

\vskip.05truein
\n{\smc Step 1.} Inequality (3.10) holds for any $C>0$ if
$$A_{n}\subset\bigl\{(x,y):\max\bigl(f_{n}(x),f_{n}(y)\bigr)\geq
     2^{n}(\log n)^2\bigr\}.$$
In this case, by (3.8),
$$\eqalign{\sum_{n}
\Pr\biggl\{\Bigl|\sum_{i,j\le 2^n}\ve_{i}\vet_{j}h_{n}(X_{i},Y_{j})\Bigr|&>
  C2^{n}\log n\biggr\}\cr
&\leq
\sum_{n}\Pr\biggl\{\exists~ i\le 2^n:\ f_{n}(X_{i})\geq 2^{n}(\log
n)^2\biggr\}\cr
&~~~~~~~~~~~+
  \sum_{n}\Pr\biggl\{\exists~ j\le 2^n:\ f_{n}(Y_{j})\geq 2^{n}(\log
n)^2\biggr\}\cr
&\le 2\sum_{n}2^{n}\Pr\bigl\{f_{n}(X)\geq 2^{n}(\log n)^2\bigr\}<\infty.\cr}$$

\vskip.05truein

\n{\smc Step 2.} Inequality (3.10) holds for any $C>0$ if
$$A_{n}\subset\bigl\{(x,y):h^{2}(x,y)\geq 2^{2n}(\log n)^2\bigr\}.$$
Indeed, by Chebyshev's inequality,
$$\eqalign{\sum_{n}
\Pr\biggl\{\Bigl|\sum_{i,j\le 2^n}\ve_{i}\vet_{j}h_{n}(X_{i},Y_{j})\Bigr|&>
  C2^{n}\log n\biggr\}\cr
&\leq
\sum_{n}{1\over C2^{n}\log n}
  E\Bigl|\sum_{i,j\le 2^n}\ve_{i}\vet_{j}h_{n}(X_{i},Y_{j})\Bigr|\cr
&=
  \sum_{n}{2^{n}\over C\log n}E|h|I_{\{|h|\geq 2^{n}\log n\}}\cr
&=
C^{-1}E|h|\sum_{n}{2^{n}\over \log n}I(|h|\geq 2^{n}\log n)\cr
&\leq
  \tilde{C}E{h^{2}\over (L_{2}|h|)^{2}}<\infty. \cr}$$

\vskip.05truein
\n{\smc Step 3}. Inequality (3.10) holds for any $C>0$ if
$$A_{n}\subset\bigl\{(x,y):2^{2n}n^{-4}\leq h^{2}(x,y)<2^{2n}(\log n)^2,
  f_{n}(x),f_{n}(y)\leq 2^{n}(\log n)^2\bigr\}.$$
If we use again Chebyshev's inequality, it suffices to prove that
$$
 \sum_{n}{E\big|\sum_{i,j\le 2^n}\ve_{i}\vet_{j}h_{n}(X_{i},Y_{j})\big|^{4}
\over
 2^{4n}(\log n)^4}<\infty.\eqno(3.11)$$
Notice however that, by iteration of Khinchin's inequality (or by direct
computation), there is $C<\infty$ (e.g. $C=18$) such that
$$\eqalign{C^{-1}E\Bigl|\sum_{i,j\le
2^n}\ve_{i}\vet_{j}h_{n}(X_{i},Y_{j})\Bigr|^{4}&\leq
 E\Bigl|\sum_{i,j\le 2^n}h^{2}_{n}(X_{i},Y_{j})\Bigr|^{2}\cr
&\le
\sum_{i,j}Eh^{4}_{n}(X_{i},Y_{j})+
  \sum_{i\neq i',j}Eh^{2}_{n}(X_{i},Y_{j})h^{2}_{n}(X_{i'},Y_{j})\cr
&
~~~~~~~~~~~~~~~~~~~~~~~~~~+\sum_{i,j\neq
j'}Eh^{2}_{n}(X_{i},Y_{j})h^{2}_{n}(X_{i},Y_{j'})\cr
&~~~~~~~~~~~~~~~~~~~~~~~~~~~+
  \sum_{i\neq i',j\neq j'}Eh^{2}_{n}(X_{i},Y_{j})
 h^{2}_{n}(X_{i'},Y_{j'}).\cr}$$
So, to prove (3.11) we have to check convergence of these four series.

First series:
$$\eqalign{\sum_{n}{2^{2n}Eh_{n}^{4}\over2^{4n}(\log n)^4}&\leq
  \sum_{n}{1\over2^{2n}(\log n)^4}Eh^{4}I_{\{h^{2}\leq 2^{2n}(\log n)^2\}}\cr
&=
Eh^{4}\sum_{n}{1\over2^{2n}(\log n)^4}I(h^{2}\leq 2^{2n}(\log n)^2)\cr
&\leq
  \tilde{C}Eh^{4}{1\over h^{2}(L_{2}|h|)^{2}}<\infty.\cr}$$

Second series: (below we use the notation
$h_{n}:=h_{n}(X,Y)$,\ $\tilde{h}_{n}=h_{n}(\tilde{X},Y)$ and
$\tilde{X}$ is an independent copy of $X$)
$$\eqalign{\sum_{n}&{2^{3n}Eh_{n}^{2}(X,Y)h_{n}^{2}(\tilde{X},Y)\over
2^{4n}\log^{4}n}=
  \sum_{n}{Eh_{n}^{2}\tilde{h}_{n}^{2}\over 2^{n}(\log n)^4}\leq
2\sum_{n}{Eh_{n}^{2}\tilde{h}_{n}^{2}I(|h|\leq |\tilde{h}|)\over
 2^{n}(\log n)^4}\cr
&\leq
 2Eh^{2}\tilde{h}^{2}I(|h|\leq |\tilde{h}|)\sum_{n}{1\over2^{n}(\log n)^4}
  I(E_{X}\min(h^{2},2^{4n})\leq 2^{n}(\log n)^2,
  \tilde{h}^{2}\leq 2^{4n})\cr
&\leq
2Eh^{2}\tilde{h}^{2}I(|h|\leq |\tilde{h}|)\sum_{n}{1\over2^{n}(\log n)^4}
  I(E_{X}\min(h^{2},\tilde{h}^{2})\leq 2^{n}(\log n)^2,|\tilde{h}|\leq
2^{2n})\cr
&  \leq
\tilde{C}Eh^{2}\tilde{h}^{2}I(|h|\leq |\tilde{h}|)
  {1\over E_{X}\min(h^{2},\tilde{h}^{2})(L_{2}|\tilde{h}|)^{2}}\leq
 \tilde{C}E{\tilde{h}^{2}\over (L_{2}|\tilde{h}|)^{2}}<\infty.\cr}$$

3rd series: convergence follows just as for the second.

4th series: here we have by (3.2)
$$\eqalign{\sum_{n}{2^{4n}(Eh_{n}^{2})^{2}\over2^{4n}(\log n)^4}&\leq
 \tilde{C} \sum_{n}{Eh_{n}^{2}\over(\log n)^3}\cr
&= \tilde{C}
Eh^{2}\sum_{n}{1\over(\log n)^3}I(2^{2n}n^{-4}\leq h^{2}(x,y)<2^{2n}
  (\log n)^2)\cr
&\leq
  \tilde{C}E{h^{2}\over(L_{2}|h|)^{2}}<\infty,\cr}$$
where we use the fact that
$$\card\bigl\{n:2^{2n}n^{-4}\leq h^{2}(x,y)<2^{2n}\log^{2}n\bigr\}\sim
2L_{2}h.$$
This completes the third Step.

\vskip.05truein
\n{\smc Step 4.} Inequality (3.10) holds for any $C>0$ if
$$A_{n}\subset \Bigl\{(x,y):h^{2}(x,y)\leq {2^{2n}\over n^{4}}, {2^{n}\over\log
n}\leq \max\bigl(f_{n}(x),f_n(y)\bigr)
                      \leq 2^{n}(\log n)^2
                      \Bigr\}.$$
We follow the proof of the previous step. The only difference is in the proof
of convergence of the fourth series. We have for $n\geq 2$
$$\eqalign{Eh_{n}^{2}&\leq 2\sum_{k=1}^{3}E\min(h^{2},2^{2n})
   I_{\{2^{n}(\log n)^{2-k}\leq f_{n}(X)\leq 2^{n}(\log n)^{3-k}\}}\cr
&\leq
\sum_{k=1}^{3}2^{n+1}(\log n)^{3-k}\Pr\bigl\{f_{n}(X)\geq
2^{n}(\log n)^{2-k}\bigr\}.\cr}$$
Thus, by (3.8),
$$\sum_{n}{Eh_{n}^{2}\over(\log n)^3}\leq \sum_{k=1}^{3}
  \sum_{n}{2^{n+1}\over(\log n)^k}\Pr\bigl\{f_{n}(X)\geq 2^{n}(\log n)^{2-k}
  \bigr\}<\infty.$$

\vskip.05truein
For the next step, we define the functions
$$g_n(x)=E_YhI_{\{|h|\ge2^nn^2\}}.\eqno(3.12)$$

\vskip.05truein

\n{\smc Step 5.} Inequality (3.1)
holds for any $C>0$ if
$$A_n\subset\bigl\{(x,y):\max\bigl(g_n(x),g_n(y)\bigr)\ge1\bigr\}.$$
Assumption (3.2) implies that $\Pr\{|h|\ge v\}\le v^{-2}L_2 v^2$. Hence,
$E|h|I_{\{|h|\ge s\}}\le\tilde{C}s^{-1}L_2s$ for $s\ge 1$. Therefore,
$$\sum_n2^n\Pr\bigl\{|g_n(X)|\ge1\bigr\}\le \tilde{C}\sum_n{Ln\over
n^2}<\infty,$$
and the same is true for $g_n(Y)$.

\vskip.05truein
\n{\smc Step 6.} Inequality (3.10) holds for any $C>0$ if
$$A_{n}\subset\Bigl\{(x,y): f_{n}(x)\geq {2^{n}\over n},\
          f_{n}(y)\geq {2^{n}\over n},\ h^{2}(x,y)\leq
          {2^{2n}\over n^{4}}\Bigr\}.$$
To see this we note first that
$$\eqalign{Eh_{n}^{2}&\leq {2^{2n}\over n^{4}}EI_{A_{n}}\leq
   {2^{2n}\over n^{4}}\Pr\Bigl\{f_{n}(X)\geq {2^{n}\over n}\Bigr\}
   \Pr\Bigl\{f_{n}(Y)\geq {2^{n}\over n}\Bigr\}\cr
&\leq
 {2^{2n}\over n^{4}}\biggl({nEf_{n}(X)\over 2^{n}}\biggr)^{2}\leq
   \tilde{C}{(\log n)^2\over n^{2}},\cr}$$
since $Ef_{n}(X)=E\min(h^{2},2^{4n})\leq \tilde{C}\log n$
by (3.2). Now we
may conclude Step 6 by Chebyshev's inequality as
$$\sum_{n}{E\big|\sum_{i,j\le 2^n}\ve_{i}\vet_{j}h_{n}(X_{i},Y_{j})\big|^2
\over
 2^{2n}(\log n)^2}\le\sum_n{Eh_n^2\over(\log n)^2}\le \tilde{C}\sum_n{1\over
n^2}<\infty.$$

\vskip.05truein
\n{\smc Step 7.} Inequality (3.10) holds for some $C>0$ if
$$A_{n}=\Bigl\{(x,y): f_{n}(x)\leq {2^{n}\over\log n},
          f_{n}(y)\leq {2^{n}\over n},g_n(x)\le 1,g_n(y)\le 1, h^{2}(x,y)\leq
          {2^{2n}\over n^{4}}\}.$$

This is the most involved step, and the only one (except for the similar Step 8
below) where we use  condition (3.3).  To prove (3.10) in this case,
we will use Prohorov's inequality (or Lemma 2.8)
together with the following four lemmas (one of which also uses Talagrands's
inequality).

\proclaim{\smc Lemma 3.3}.
For all $n\in\bf N$,
$$\Pr\biggl\{\Bigl|\sum_{i\le 2^n}\ve_{i}h_{n}(X_{i},Y)\Bigr|
  \geq 2^{n+4}\biggr\}\leq 2^{-4n}$$
and
$$\sum_n\Pr\biggl\{\max_{1\le j\le 2^n}\Bigl|\sum_{i\le
2^n}\ve_{i}h_{n}(X_{i},Y_j)\Bigr|
  \geq 2^{n+4}\biggr\}<\infty$$

\n{\smc Proof.} We note that $A_{n}\subset \bigl\{(x,y):|h(x,y)|\leq
n^{-1}2^{n},  f_{n}(y)\leq n^{-1}2^{n}\bigr\}$ and then apply Bernstein's
inequality  or  Prohorov's inequality
to obtain that,
for any $Y$,
$${\Pr}_{X}\biggl\{\Bigl|\sum_{i\leq 2^{n}}\ve_{i}h_{n}(X_{i},Y)\Bigr|\geq
2^{n+4}
  \biggr\}\leq e^{-4n},$$
which clearly implies the Lemma. (Lemma 2.8 instead of Bernstein's or
Prohorov's inequality would simply change multiplicative constants.)
\qed

\vskip.05truein
Before formulating the next lemma it is convenient to define a sequence $c_{n}$
by the formula
$$c_{n}=Eh^{2}I_{\{2^{n}n^{-2}< |h|\leq 2^{n}n^{2}\}},\ \ n\in\bf
N.\eqno(3.13)$$

\n{\smc Lemma 3.4}.
We have
$$\sum_{n}\exp\biggl(-{2\log n\over\sqrt{1+c_{n}}}\biggr)<\infty.$$

\n{\smc Proof.}
Condition (3.2) implies that, for any $k\geq 2$,
$$\sum_{k\leq \log n\leq k+1}c_{n}\leq \tilde{C}kE|h|^{2}
  I_{\{|h|\leq 2^{e^{k+1}}(e^{k+1})^{2}\}}\leq \tilde{C}k^{2},$$
(where the second constant is different from the first) since the largest
number
of intervals
$I_n=[n^{-2}2^n,n^22^n]$, $k\le \log n\le k+1$, that can overlap with any
given
one of them is not larger than
$6(k+1)$.
Hence,
$$\card\{n: k\leq \log n\leq k+1, c_{n}\geq 1\}\leq \tilde{C}k^{2}.$$
Condition (3.2) also implies
$c_{n}\leq 2\log n$ (note that $c_1=0$). So,
$$\eqalign{\sum_{n}\exp\biggl(-{2\log n\over\sqrt{1+c_{n}}}\biggr)&\leq
  \sum_{n}\exp\bigl(-\sqrt{2}\log n\bigr)+
  \sum_{c_{n}\geq 1}\exp\biggl(-{2\log n\over\sqrt{1+2\log n}}\biggr)\cr
&\leq
\sum_{n}\exp\bigl(-\sqrt{2}\log
n\bigr)+\sum_{k}\tilde{C}k^{2}\exp\bigl(-\sqrt{k}\bigr)<\infty.\cr}$$
\rightline{\qed}

\vskip.05truein
The following lemma is well known but a proof is provided for the reader's
convenience.

\proclaim{\smc Lemma 3.5}. If a kernel $k$ satisfies
$E_X|k(X,y)|\le 1$ and $E_Y|k(x,Y)|\le 1$ a.s., then $k$ defines an operator on
$L_2({\cal L}(X))$ with norm bounded by 1, that is, condition (3.3) holds for
$h=k$ and
$C=1$ (and therefore so does condition (1.5)).

\n{\smc Proof}. We need to
check that
$$|E_XE_Y k(X,Y)f(X)g(Y)|\le [Ef^2(X)Eg^2(Y)]^{1/2}
$$
whenever $\|f\|_\infty,$ $\|g\|_\infty<\infty$.
But, assuming (without loss of generality) that $k$, $f$ and $g$ are
nonnegative,
$$\eqalign{E_XE_Y k(X,Y)f(X)g(Y)&=E_X \Bigl[f(X)E_Y\bigl(k^{1/2}(X,Y)
k^{1/2}(X,Y)g(Y)\bigr)\Bigr]\cr
&\le E_X \Bigl[f(X) (E_Yk(X,Y))^{1/2} (E_Y k(X,Y)g^2(Y))^{1/2}\Bigr]\cr
&\le E_X
\Bigl[f(X)\bigl(E_Y k(X,Y)g^2(Y)\bigr)^{1/2}\bigr]\cr
&\le \bigl(E_X f^2(X)\bigr)^{1/2} \Bigl[E_X\bigl(E_Y
k(X,Y)g^2(Y)\bigr)\Bigr]^{1/2}.\cr}$$
and now the inequality follows by applying Fubini and using $E_Xk(X,Y)\le 1$.
\qed

\proclaim{\smc Lemma 3.6}.  There exists $C_{1}<\infty$ such that
$$\sum_{n}\Pr\biggl\{E_{Y}\Bigl(\sum_{i\le 2^n}\ve_{i}h_{n}(X_{i},Y)
   \Bigr)^{2}\geq C_{1}\sqrt{1+c_{n}}~2^{n}\log n\biggr\}<\infty.$$

\n{\smc Proof.}  Let
$H_Y$ be
$L_2(\Omega,
\sigma(Y),\Pr)$, that is,
$H_Y$ is the space of all square integable random variables $f(Y)$  where
$f$ is
a Borel measurable function.
Let ${\bf X}_i:=\ve_ih_n(X_i,Y)$ for $i=1,\dots,2^n$. Then, ${\bf X}_i$ are
symmetric i.i.d. random vectors with values in $H_Y$.
We define
$$Z=\sup_{f\in{\cal F}}\sum_{i=1}^{2^n}f({\bf X}_i)
=\biggl[E_{Y}\Bigl(\sum_{i\le 2^n}\ve_{i}h_{n}(X_{i},Y)
   \Bigr)^{2}\biggr]^{1/2},$$
where $\cal F$ is a countable dense subset of the unit ball of $H'_Y=H_Y$ and
we write $f(\cdot):=\langle f,\cdot\rangle$. We will apply Lemma 2.8 to $Z$.
For this, we must estimate $EZ$ and determine suitable $U$ and $\sigma^2$. We
have
$$EZ\le\bigl(EZ^2\bigr)^{1/2}=\bigl[2^nEh_n^2\bigr]^{1/2}\le
\sqrt{2^n\log n}\eqno(3.14)$$
by (3.2). Since
$$\sup_{f\in{\cal F}}\big|f({\bf X}_i(\omega))\big|=\|{\bf
X}_i(\omega)\|_Y
=\sqrt{E_Yh_n^2(X_i(\omega),Y)}\le \sqrt{2^n\over \log n},$$
we can take
$$U=\sqrt{2^n\over \log n}\eqno(3.15)$$
in Lemma 2.8 for $Z$. Moreover, for each $f\in{\cal F}$,
$$Ef^2({\bf X}_i)=E\bigl(E_Yh_n(X_i,Y)f(Y)\bigr)^2\le 3\sum_{i=1}^3
E\bigl(E_Yh_n^{(i)}(X_i,Y)f(Y)\bigr)^2,$$
where
$$h_n^{(1)}:=hI_{B_n},\ h_n^{(2)}:=hI_{B_n\cap\{2^nn^{-2}<|h|\le2^nn^2\}},\
h_n^{(3)}:=hI_{B_n\cap\{|h|\ge2^nn^2\}},$$
with
$$B_n:=\Bigl\{(x,y):f_n(x)\le{2^n\over \log n},f_n(y)\le {2^n\over n},
g_n(x)\le 1, g_n(y)\le 1\Bigr\}, $$
since
$$h_n=h_n^{(1)}-h_n^{(2)}-h_n^{(3)}.$$
Now,
$$E\bigl(E_Yh_n^{(1)}(X_i,Y)f(Y)\bigr)^2\le 1$$
by condition (1.5) (which is equivalent to (1.3)=(3.3)),
$$E\bigl(E_Yh_n^{(2)}(X_i,Y)f(Y)\bigr)^2\le E(h_n^{(2)})^2\le c_n$$
by Cauchy-Schwartz and the definition of $c_n$ in (3.13), and
$$E\bigl(E_Yh_n^{(3)}(X_i,Y)f(Y)\bigr)^2\le 1$$
by Lemma 3.5 (see (1.5) once more). Therefore, we can take $\sigma^2$ in Lemma
2.8 for $Z$ to be
$$\sigma^2=3\cdot2^n(2+c_n)<6\cdot2^n(1+c_n).\eqno(3.16)$$
Then, on account of (3.14)-(3.16),
Lemma 2.8 gives, with $C_2=(\sqrt{C_1}-1)^2$,
$$\eqalign{\Pr\biggl\{
E_{Y}\Bigl|\sum_{i\le2^n}\ve_{i}h_{n}&(X_{i},Y)
   \Bigr|^{2}\geq C_{1}\sqrt{1+c_{n}}~2^{n}\log n\biggr\}\cr
&=\Pr\biggl\{Z\ge\sqrt{C_{1}\sqrt{1+c_{n}}~2^{n}\log n}\biggr\}\cr
&\le \Pr\biggl\{Z-EZ\ge\sqrt{C_{2}\sqrt{1+c_{n}}~2^{n}\log n}\biggr\}\cr
&\le  K\exp\biggl(-{\sqrt{C_2}\over K}{\root 4 \of{1+c_{n}}}~\log
n~\log\Bigl(1+{\sqrt{C_2}{\root 4 \of{1+c_{n}}}\over
6(1+c_n)+8}\Bigr)\biggr)\cr
&\le K\exp\biggl(-{\sqrt{C_2}\over K}{\root 4 \of{1+c_{n}}}~\log
n~\log\Bigl(1+{\sqrt{C_2}\over
14(1+c_n)^{3/4}}\Bigr)\biggr)\cr
&\le K\exp\biggl(-{\sqrt{C_2}\over K}\log\Bigl(1+{\sqrt{C_2}\over
14}\Bigr){\log n\over\sqrt{1+c_n}}\biggr),\cr}$$
where in the last line we have used that the function $x^{-1}\log(1+x)$ is
monotone decreasing. Taking $K^{-1}\sqrt{C_2}\log(1+\sqrt{C_2}/14)\ge 2$ yields
the bound
$$\Pr\biggl\{
E_{Y}\Bigl|\sum_{i\le 2^n}\ve_{i}h_{n}(X_{i},Y)
   \Bigr|^{2}\geq C_{1}\sqrt{1+c_{n}}~2^{n}\log
n\biggr\}\le K\exp\biggl(-{2\log n\over\sqrt{1+c_n}}\biggr)$$
and Lemma 3.6 follows from Lemma 3.4. \qed

\vskip.05truein

Now we complete the proof of Step 7. For $n$ fixed, set
$$d(y):=\sum_{i\le 2^n}\varepsilon_ih_n(X_i,y)\ \ {\rm and}\ \
\tilde{d_j}:=\vet_jd(Y_j)I_{\{|d|\le 2^{n+4},E_Yd^2(Y)\le C_12^n(\log
n)\sqrt{1+c_n}\}}(Y_j)$$
for $1\le j\le 2^n$. Then,
$$\eqalign{\Pr\biggl\{\Big|\sum_{i,j\le
2^n}\ve_i\vet_jh_n&(X_i,Y_j)\Big|>C2^n\log n\biggr\}
=\Pr\biggl\{\Big|\sum_{j\le
2^n}\vet_jd(Y_j)\Big|>C2^n\log n\biggr\}\cr
&\le\Pr\bigl\{\exists~j\le 2^n: \tilde{d}_j\ne d(Y_j)\bigr\}
+\Pr\biggl\{\Big|\sum_{j\le
2^n}\vet_j\tilde{d}_j\Big|>C2^n\log n\biggr\}\cr
&:=I_n+II_n.\cr}$$
But,
$$\eqalign{I_n\le\Pr&\biggl\{\max_{j\le 2^n}\Big|\sum_{i\le
2^n}\ve_ih_n(X_i,Y_j)\Big|>2^{n+4}\biggr\}\cr
&+\Pr\biggl\{E_Y\Big|
\sum_{i\le
2^n}\ve_ih_n(X_i,Y)\Big|^2>C_12^n(\log
n)\sqrt{1+c_n}\biggr\}\cr}$$
and Lemma 3.3 and Lemma 3.6 show that
$$\sum_nI_n<\infty.\eqno(3.17)$$
To estimate $II_n$ we can apply Bernstein's or Prokhorov's inequality
conditionally on the sequence $\{X_i\}$. For
convenience we will use Lemma 2.8. We can take
$U=2^{n+4}$ and
$V=C_12^{2n}(\log n)\sqrt{1+c_n}$ to get
$$\eqalign{{\Pr}_Y\biggl\{\Big|\sum_{j\le
2^n}\vet_j\tilde{d}_j\Big|&>C2^n\log n\biggr\}\cr
&\le K\exp\biggl(-{1\over
K}{C2^n\log n\over 2^{n+4}}\log\Bigl(1+{C2^{2n+4}\log n \over C_12^{2n}(\log
n)\sqrt{1+c_n}}\Bigr)\biggr)\cr
&\le K\exp\biggl(-{C\over 2^4K}\log\Bigl(1+{2^4C\over C_1}\Bigr){\log
n\over\sqrt{1+c_n}}\Bigr)\biggr).\cr}$$
Taking $C$ so that
$${C\over 2^4K}\log\Bigl(1+{2^4C\over C_1}\Bigr)\ge 2$$
shows, by Lemma 3.4, that
$$\sum_nII_n<\infty.\eqno(3.18)$$
(3.17) and (3.18) complete the proof of Step 7.

\vskip.05truein
\n{\smc Step 8.} Inequality (3.10) holds for some $C<\infty$ if
$$A_{n}=\Bigl\{(x,y): f_{n}(x)\leq {2^{n}\over n},f_{n}(y)\leq {2^{n}\over\log
n}, g_n(x)\le 1,g_n(y)\le 1, h^{2}(x,y)\leq
          {2^{2n}\over n^{4}}\}.$$
This can be done in the same way  as Step 7.

\vskip.05truein
It is clear that we can write $S\times S=\cup_{i=1}^8A_n^i$ with $A_n^1,\dots,
A_n^8$ disjoint, and $A_n^i$ satisfying the conditions in Step i for each $n$.
Then, $h=\sum_{i=1}^8hI_{A_n^i}=\sum_{i=1}^8h_n^i$. Since for each $i$ the
kernels $h_n^i$ satisfy condition (3.10) for some $C<\infty$, it follows by the
triangle inequality that  the series (3.4) for $h$ converges for some
$C<\infty$, proving the sufficiency part of Theorem 3.1.
\qed

\vskip.1truein
\n{\bf 4. Canonical kernels.} In this section we show that, for canonical
kernels,
the LIL (1.1) is equivalent to the decoupled and randomized LIL. The
preliminary results in Section 2(B) yield that the regular LIL implies the
decoupled and randomized one.  The converse implication, however, seems to
require Theorem 3.1. The first step consists of the following simple
inequality,
rooted in known symmetrization techniques.

\n\proclaim{\smc Lemma 4.1}. For any kernel $h$, and for any $n\in\bf N$ and
$t>0$, we have
$$\eqalign{\Pr\biggl\{\Bigl|\sum_{i,j\le n}h(X_{i},Y_{j})\Bigr|\geq 10t\biggr\}
\leq
    16\Pr&\biggl\{\Bigl|\sum_{i,j\le n}\ve_{i}\vet_{j}h(X_{i},Y_{j})\Bigr|\geq
t\biggr\}\cr
&+ 4\Pr\biggl\{E_{Y}\Bigl|\sum_{i,j\le n}\ve_{i}h(X_{i},Y_{j})\Bigr|\geq
t\biggr\}\cr
&~~~~~~~~+
 \Pr\biggl\{E_{X}\Bigl|\sum_{i,j\le n}h(X_{i},Y_{j})\Bigr|\geq
t\biggr\}.\cr}$$

\n{\smc Proof.} Let $\{Z_{i}\}$ be a sequence of independent random variables
such that $E|\sum_{i}Z_{i}|\leq s$ and let $\{Z_{i}^{'}\}$ be an
independent copy
of $\{Z_{i}\}$. Then, by Chebyshev's inequality,
$\Pr\bigl\{|\sum_{i}Z_{i}^{'}|\leq 2s\bigr\}\geq 1/2$. So, for any $t>0$,
$$\eqalign{\Pr\Bigl\{\Bigl|\sum_{i}Z_{i}\Bigr|\geq
2t+2s\biggr\}&\le2\Pr\biggl\{
\Big|\sum_iZ'_i\Big|\le 2s,\Big|\sum_iZ_i\Big|\ge 2t+2s\biggr\}\cr &\leq
  2\Pr\biggl\{\Bigl|\sum_{i}(Z_{i}-Z_{i}^{'})\Bigr|\geq 2t\biggr\}\cr
&=
  2\Pr\biggl\{\Bigl|\sum_{i}\ve_{i}(Z_{i}-Z_{i}^{'})\Bigr|\geq 2t\biggr\}\cr
&\leq
2\Pr\biggl\{\Bigl|\sum_{i}\ve_{i}Z_{i}\Bigr|\geq t\biggr\}+
  2\Pr\biggl\{\Bigl|\sum_{i}\ve_{i}Z_{i}^{'}\Bigr|\geq t\biggr\}\cr
&=
  4\Pr\biggl\{\Bigl|\sum_{i}\ve_{i}Z_{i}\Bigr|\geq t\biggr\}.\cr}$$
Using the above inequality conditionally we get
$$\Pr\biggl\{\Bigl|\sum_{i,j}h(X_{i},Y_{j})\Bigr|\geq 10t\biggr\}\leq
4\Pr\biggl\{\Bigl|\sum_{i,j}\ve_{i}h(X_{i},Y_{j})\Bigr|\geq 4t\biggr\}+
  \Pr\biggl\{E_{X}\Bigl|\sum_{i,j}h(X_{i},Y_{j})\Bigr|\geq t\biggr\}$$
and
$$\eqalign{\Pr\biggl\{\Bigl|\sum_{i,j}\ve_{i}h(X_{i},Y_{j})\Bigr|\geq
4t\biggr\}\leq
4\Pr\biggl\{\Bigl|&\sum_{i,j}\ve_{i}\vet_{j}h(X_{i},Y_{j})\Bigr|\geq
t\biggr\}\cr &+
  \Pr\biggl\{E_{Y}\Bigl|\sum_{i,j}\ve_{i}h(X_{i},Y_{j})\Bigr|\geq
t\biggr\}.\cr}$$
\rightline{\qed}

\vskip.05truein
The next lemma shows that if the second moment and the conditional second
moment of a canonical kernel $h$ are suitably truncated, then
Talagrand's inequality (Lemma 2.8) allows control of the last two terms on the
right hand side of the inequality in Lemma 4.1.

\proclaim{\smc Lemma 4.2}.
 Let $h$ be a canonical kernel such that
$$ Eh^{2}(X,Y)\leq c^2\log n \ \  and \ \ E_{Y}h^{2}(X,Y)\leq c^22^{n}
   \ \ \ X-a.s.$$
for some $c<\infty$.
 Then we have that, for some universal constant $C$,
 $$\Pr\biggl\{E_{Y}\Bigl|\sum_{i,j\leq 2^{n}}h(X_{i},Y_{j})\Bigr|\geq
  cC2^{n}\log n\biggr\}\leq n^{-2}.$$

\n{\smc Proof.} We can assume $c=1$.
If we define
$$Z:=E_{Y}\Bigl|\sum_{i,j\leq 2^{n}}h(X_{i},Y_{j})\Bigr|$$
then
$$Z=\sup\Bigl\{
  \sum_{i\leq 2^{n}}E_{Y}\bigl(\sum_{j\leq 2^{n}}h(X_{i},Y_{j})g({\bf Y})\Bigr)
  \Bigr\},$$
where the supremum is taken over all $g({\bf Y})=g(Y_{1},\ldots,Y_{2^{n}})$
with
$\|g\|_{\infty}\leq 1$, actually over a countable $L_1$-norm determining
subset of such functions.
Thus
$Z$ has the same form as in Lemma 2.8. Then, since
$$
 \Big\|E_{Y}\big|\sum_{j=1}^{2^n}h(x,Y_{j})\big|\Big\|_{\infty}
 \leq \Big\|\Bigl(\sum_{j=1}^{2^n}E_{Y}h^2(x,Y_{j})\Bigr)^{1/2}\Big\|_{\infty}
 \leq 2^n$$
and
$$
\sum_{i=1}^{2^n}E\Bigl(E_Y\big|\sum_{j=1}^{2^n}h(X_i,Y_j)\big|\Bigr)
  ^{2}
\leq
  \sum_{i=1}^{2^n}E\Bigl(\sum_{j=1}^{2^n}E_Yh^{2}(X_i,Y_j)\Bigr)
=2^{2n}Eh^{2}\leq 2^{2n}\log n,$$
we can take
$$U=2^n\ \ {\rm and}\ \ V=2^{2n}\log n\eqno(4.1)$$
in Talagrand's exponential bound for $Z$.
Moreover
$$EZ\leq \Bigl(E\bigl|\sum_{i,j\leq 2^{n}}h(X_{i},Y_{j})\bigr|^{2}\Bigr)^{1/2}
  =2^{n}\bigl(Eh^{2}\bigr)^{1/2}\leq 2^{n}\log n.\eqno(4.2)$$
Now the statement follows by (4.1), (4.2) and the exponential bound in
Lemma 2.8. \qed

\vskip.05truein
The following lemma will allow us to carry out truncations for canonical
kernels exactly in the same way as we did for randomized kernels in the first
four steps of the sufficiency proof of Theorem 3.1.

\proclaim{\smc Lemma 4.3}.
For any integrable kernel $h$, $n\in\bf N$ and $p\geq 1$ we have
$$\Bigl\|\sum_{i,j\le n}\pi_2h(X_{i},Y_{j})\Bigl\|_{p}\leq 4
  \Bigr\|\sum_{i,j\le n}\ve_{i}\vet_{j}h(X_{i},Y_{j})\Bigr\|_{p}.$$

\n{\smc Proof.} Since $\pi_2h$ is canonical, by Jensen's inequality we have
that,
 for all $\{Y_{j}\}$,
$$\eqalign{ E_{X}\Bigl|\sum_{i,j\le 2^n}\pi_2h(X_{i},Y_{j})\Bigr|^{p}&\leq
   E_{X}\Bigl|\sum_{i,j\le 2^n}\bigl(\pi_2h(X_{i},Y_{j})-
   \pi_2h(X_{i}^{'},Y_{j})\bigr)\Bigr|^{p}\cr
&=
E_{X}\Bigl|\sum_{i,j\le 2^n}\ve_{i}\bigl(\pi_2h(X_{i},Y_{j})-
   \pi_2h(X_{i}^{'},Y_{j})\bigr)\Bigr|^{p}\cr
&=
 E_{X}\Bigl|\sum_{i,j\le 2^n}\ve_{i}\bigl(h(X_{i},Y_{j})-E_{Y}h(X_{i},Y_{j})\cr
&~~~~~~~~~~~~~~~~~~~~~~~~~~~~~~~-
    h(X_{i}^{'},Y_{j})+E_{Y}h(X_{i}^{'},Y_{j})\bigr)\Bigr|^{p}.\cr}$$
 Thus, by the triangle inequality,
$$\eqalign{\Bigl\|\sum_{i,j\le 2^n}\pi_2h(X_{i},Y_{j})\Bigr\|_{p}&\leq
 \Bigl\|\sum_{i,j\le 2^n}\ve_{i}\bigl(h(X_{i},Y_{j})-E_{Y}h(X_{i},Y_{j})\bigr)
   \Bigr\|_{p}\cr
&~~~~~~~~~~+
   \Bigl\|\sum_{i,j\le
2^n}\ve_{i}\bigl(h(X_{i}^{'},Y_{j})-E_{Y}h(X_{i}^{'},Y_{j})
   \bigr)\Bigr\|_{p}\cr
&=
 2\Bigl\|\sum_{i,j\le 2^n}\ve_{i}\bigl(h(X_{i},Y_{j})-E_{Y}h(X_{i},Y_{j})\bigr)
   \Bigr\|_{p}.\cr}$$
 In a similar way we may prove that
$$\Big\|\sum_{i,j\le 2^n}\ve_{i}\bigl(h(X_{i},Y_{j})-E_{Y}h(X_{i},Y_{j})\bigr)
   \Bigr\|_{p}\leq
   2\Bigl\|\sum_{i,j\le 2^n}\ve_{i}\vet_{j}h(X_{i},Y_{j})\Bigr\|_{p}.$$
\rightline{\qed}

\vskip.05truein
Now we can prove the main result of this section.

\proclaim{\smc Theorem 4.4}.
For any canonical kernel $h$ the following two conditions are
equivalent:
$$
 \limsup_{n\rightarrow\infty}{1\over n\log\log n}\Bigl|
 \sum_{1\leq i< j\leq n}h(X_i,X_j)\Bigr|<\infty \ \ a.s.\eqno(4.3)
$$
and
$$
 \limsup_{n\rightarrow\infty}{1\over n\log\log n}\Bigl|
 \sum_{1\leq i,j\leq n}\ve_{i}\vet_{j}h(X_i,Y_j)\Bigr|<\infty \ \ a.s.
\eqno(4.4)$$

Here, again, each of the two limsups is a.s. bounded by a universal constant
times the other.

\n{\smc Proof.}
(4.3) implies (4.4) (even without degeneracy of the kernel)
by Lemma 2.1(b).

To prove the opposite implication, by Corollary 2.6 it is enough
to show that if (4.4) holds (which is equivalent to the two conditions
(3.2) and
(3.3) by Theorem 3.1),  then
$$
  \sum_{n}\Pr\biggl\{\Bigl|\sum_{i,j\leq 2^{n}}h(X_{i},Y_{j})\Bigr|\geq
  C2^{n}\log n\biggr\}<\infty.
$$
Since $h$ is canonical, we may replace $h$ by
$\pi_2h$ in this series ($h=\pi_2h$). As in the case
of decoupled and randomized kernels, convergence of the series will follow
in a few steps by showing that
$$
 \sum_{n}\Pr\biggl\{\Bigl|\sum_{i,j\leq 2^{n}}\pi_2h_{n}(X_{i},Y_{j})\Bigr|
 \geq C2^{n}\log n\biggr)<\infty,\eqno(4.5)
$$
where $h_{n}=hI_{A_{n}}$ for suitably chosen sequences of sets $A_{n}$. We can
assume, as in Theorem 3.1, that  $C=1$ in conditions (3.2) and (3.3).

\vskip.05truein

\n{\smc Step 1.} The series in (4.5) converges for $$A_{n}=\bigl\{(x,y):\
f_{n}(x)>2^{n}(\log n)^2
\ \ {\rm or }\ \
              f_{n}(y)>2^{n}(\log n)^2\bigr\}.$$
By the degeneracy of $h$ we have
$$\eqalignno{|Eh_{n}|&=\big|EhI_{\{f_n(x)>2^{n}(\log
n)^2\}}+EhI_{\{f_n(y)>2^{n}(\log n)^2\}}&\cr
&~~~~~~~~~~~~~~~~~~~~~~~~~~~~~~~~~~- EhI_{\{f_{n}(x)>2^{n}(\log
n)^2,f_{n}(y)>2^{n}(\log n)^2\}}\big|&\cr
&=\bigl|EhI_{\{f_{n}(x)>2^{n}(\log
n)^2,f_{n}(y)>2^{n}(\log n)^2\}}
  \bigr|&\cr
&\leq
\Pr\bigl\{f_{n}(X)>2^{n}(\log n)^2\bigr\}^{1/2}
 \Pr\bigl\{f_{n}(Y)>2^{n}(\log
n)^2\bigr\}^{1/2}&\cr
&\leq\tilde{C}2^{-n},&(4.6)\cr}$$
where the last two inequalities follow by (3.3) and (3.8)
respectively.
We also have
$$\pi_1h_{n}(x)=\pi_1hI_{\{f_{n}(y)>2^{n}(\log n)^2,
   f_{n}(x)\leq 2^{n}(\log n)^2\}}(x),$$
as can be seen using the decomposition of $h_n$ given in the first line
of (4.6) together with the fact that
$E_YhI_{\{f_n(x)>2^{n}(\log n)^2\}}=0$. Thus, by Chebyshev's inequalty,
$$\eqalignno{\sum_{n}\Pr\biggl\{\Bigl|\sum_{i\leq
2^{n}}&\pi_1h_{n}(X_{i})\Bigr|\geq
  c\log n\biggr\}&\cr
&\leq
\sum_{n}{2^{n}\over c^{2}(\log n)^2}E\Bigl|\pi_1h
  I_{\{f_{n}(y)>2^{n}(\log n)^2, f_{n}(x)\le2^n(\log n)^2\}}(X)\Big|^2&\cr
&\le
\sum_{n}{2^{n}\over c^{2}(\log n)^2}E_{X}\Bigl[\bigl(E_{Y}
  hI_{\{f_{n}(y)>2^{n}(\log n)^2\}}\bigr)^{2}I_{\{f_n(X)\le2^{n}(\log
n)^2\}}\Bigr]&\cr
&\le
\sum_{n}{2^{n}\over c^{2}(\log n)^2}E_{X}\bigl(E_{Y}
  hI_{\{f_{n}(y)>2^{n}(\log n)^2\}}\bigr)^{2}&\cr
&\leq
  \sum_{n}{2^{n}\over c^2(\log
n)^2}\Pr\bigl\{f_{n}(Y)>2^{n}(\log n)^2\bigr\}<
  \infty,&(4.7)\cr}$$
where in the last line we used (1.5) with $C=1$ (that is, condition (3.3)) and
(3.8). Finally, as in step 1 of the proof of sufficiency of the symmetrized
LIL,
$$
 \sum_{n}\Pr\biggl\{\Bigl|\sum_{i,j\leq 2^{n}}h_{n}(X_{i},Y_{j})\Bigr|\geq
 C2^{n}\log n\biggr\}<\infty.\eqno(4.8)
$$
Inequalities (4.6)-(4.8) imply (4.5) by Hoeffding's
decomposition ((2.1)).

\vskip.05truein
\n{\smc Step 2.} The series in (4.5) converges for
$$\eqalign{A_{n}\subset\bigl\{(x,y): |h(x,y)|>2^{n}\log n
\ \ {\rm or }&\ \
              f_{n}(x)>2^{n} \ \ {\rm or }\ \ f_{n}(y)>2^{n}\bigr\}\cr
&\cap
\bigl\{(x,y): \max(f_n(x),f_n(y))\le 2^n(\log n)^2\bigr\}.\cr}$$
To prove this we may proceed just as in steps 2-4 of the proof of the
symmetrized  LIL, with only formal changes: note that in steps 2-4 there we
used
only Chebyshev's inequality to bound probabilities; thus Lemma 4.3
reduces proving inequality (4.5) here to steps 2-4 in that proof, where the
lower bounds for $h$ and $f_n$ are even smaller.

\vskip.05truein
\n{\smc Step 3.} The series in (4.5) converges for
$$A_{n}=\{(x,y):|h(x,y)|\leq
2^{n}\log n,
f_{n}(x)\le 2^n,f_{n}(y)\leq 2^{n}\}.$$
The LIL (4.4) implies that
$$
  \sum_{n}\Pr\biggl\{\Bigl|\sum_{i,j\leq 2^{n}}\ve_{i}\vet_{j}h(X_{i},Y_{j})
  \Bigr|\geq C2^{n}\log n\biggr\}<\infty
$$
for some $C<\infty$ by Lemma 2.2(c). Steps 1-4
from the proof of sufficiency in Theorem 3.1 show that
$$\sum_{n}\Pr\biggl\{\Bigl|\sum_{i,j\leq
2^{n}}\ve_{i}\vet_{j}hI_{D_n}(X_{i},Y_{j})
  \Bigr|\geq C2^{n}\log n\biggr\}<\infty,$$
for any $D_n\subset\bigl\{(x,y):|h(x,y)|>2^n/n\ \ {\rm  or}\ \
\max(f_n(x),f_n(y))>2^n/\log n\bigr\}$, in particular for $D_n=A_n^c$.
Therefore
we have
$$
\sum_{n}\Pr\biggl\{\Bigl|\sum_{i,j\leq 2^{n}}\ve_{i}\vet_{j}h_{n}(X_{i},Y_{j})
  \Bigr|\geq C2^{n}\log n\biggr\}<\infty\eqno(4.9)
$$
for some $C<\infty$. In order to deduce (4.5) from (4.9) we show first
that we can replace $h_n$  by
$\pi_2h_n$ in (4.9), and then apply Lemmas 4.1 and 4.2 to $\pi_2h_n$. So, we
begin by proving (4.9) for $h_n-\pi_2h_n$ or, what is the same by Hoeffding's
decomposition, we prove (4.9) with
$h_n$ replaced by
$\pi_1h_n$  and by
$Eh_n$.  We can write
$h_n$ as
$$\eqalign{h_{n}=h-hI_{\{f_{n}(x)>2^{n}\}}&-hI_{\{f_{n}(y)>2^{n}\}}\cr
  &+hI_{\{f_{n}(x)>2^{n},f_{n}(y)>2^{n}\}}-
hI_{\{|h|>2^{n}\log n, f_{n}(x)\le2^n,f_{n}(y)\leq 2^{n}\}}.\cr}$$
Then,  by the degeneracy of $h$ and (3.3) we have
$$\eqalign{\Bigl|\sum_{i,j\leq 2^{n}}\ve_{i}\vet_{j}&Eh_{n}\Bigr|\leq
2^{2n}\Bigl(
   \bigl|EhI_{\{f_{n}(x)>2^{n},f_{n}(y)>2^{n}\}}\bigr|+
   E|h|I_{\{|h|>2^{n}\log n\}}\Bigr)\cr
&\leq
2^{2n}\Bigl(\Pr\bigl\{f_{n}(X)>2^{n}\bigr\}^{1/2}
\Pr\bigl\{f_{n}(Y)>2^{n}\bigr\}^{1/2}+E|h|I_{\{|h|>2^{n}\log
n\}}\Bigr).\cr}$$
Now, we note that (3.2) implies $E|h|I_{\{|h|>2^{n}\log n\}}\le\tilde{C}2^{-n}$
(as $\Pr\{|h|>u\}\le u^{-2}L_2u$) and
$$\Pr\bigl\{f_{n}(X)>2^{n}\bigr\}\le{E(h^2\wedge2^{4n})\over2^n}\le
\tilde{C}{Ln\over
2^{n}}.$$
Hence,
$$\Bigl|\sum_{i,j\leq 2^{n}}\ve_{i}\vet_{j}Eh_{n}\Bigr|\le
\tilde{C}2^{n}\log n.
\eqno(4.10)$$
The above decomposition of $h_n$ together with the degeneracy of $h$ also give
$$\pi_1h_{n}(x)=-\pi_1hI_{\{f_{n}(y)>2^{n}, f_{n}(x)\leq 2^{n}\}}(x)-
  \pi_1hI_{\{|h|>2^{n}\log n, f_{n}(x),f_{n}(y)\leq 2^{n}\}}(x).$$
So, by Chebyshev's inequality and (3.2), we have
$$\eqalignno{\sum_{n}\Pr\biggl\{\Bigl|\sum_{i,j\leq 2^{n}}\ve_{i}\vet_{j}\pi_1
  h&I_{\{|h(x,y)|>2^{n}\log n, f_{n}(x),f_{n}(y)\leq 2^{n}\}}(X_{i})\Bigr|
  \geq c2^{n}\log n\biggr\}&\cr
&\leq
\sum_{n}{2^{n}\over c\log n}E\bigl|\pi_1
  hI_{\{|h(x,y)|>2^{n}\log n, f_{n}(x),f_{n}(y)\leq 2^{n}\}}(X)\bigr|&\cr
&\leq
\sum_{n}{1\over c\log n}2^{n+1}E|h|I_{\{|h|>2^{n}\log n\}}&\cr
&\leq
 c^{-1}E|h|\sum_{n}{2^{n+1}\over\log n}I(|h|>2^{n}\log n)&\cr
&\leq
\tilde{C}E{h^{2}\over(L_{2}|h|)^{2}}<\infty.&(4.11)\cr}$$
Also, by Chebyshev's inequality, (1.5) with $C=1$ and (3.8),
$$\eqalignno{\sum_{n}\Pr\biggl\{\Bigl|\sum_{i,j\leq 2^{n}}\ve_{i}\vet_{j}\pi_1
  h&I_{\{f_{n}(y)>2^{n}, f_{n}(x)\leq 2^{n}\}}(X_{i})\Bigr|
  \geq c2^{n}\log n\biggr\}&\cr
&\leq
\sum_{n}{1\over c^{2}\log^{2}n}E\Bigl|\pi_1
  hI_{\{f_{n}(y)>2^{n}, f_{n}(x)\leq 2^{n}\}}(X)\Bigr|^{2}&\cr
&\leq
\sum_{n}{1\over
c^{2}\log^{2}n}E_{X}\Bigl(E_{Y}hI(f_{n}(y)>2^{n})\Bigr)^{2}&\cr
&\leq
 \sum_{n}{1\over
c^{2}\log^{2}n}\Pr\bigl\{f_{n}(Y)>2^{n}\bigr\}<\infty.&(4.12)\cr}$$
Inequalities (4.9)-(4.12) imply, by the Hoeffding's
decomposition,
$$
  \sum_{n}\Pr\biggl\{\Bigl|\sum_{i,j\leq 2^{n}}\ve_{i}\vet_{j}\pi_2
  h_{n}(X_{i},Y_{j})\Bigr|\geq C2^{n}\log n\biggr\}<\infty\eqno(4.13)
$$
for some $C<\infty$.
By (3.2), $E(\pi_2h_n)^2\le Eh_n^2\le \tilde{C}\log n$,   and, by
the definition of $A_n$ and (3.2),
$E_Y(\pi_2h_n)^2(x)\le 2E_Yh_n^2+2Eh_n^2\le 2^{n+1}+\tilde{C}\log n$, and
likewise for
$E_X(\pi_2h_n)^2$. Then, it follows from Lemma 4.2 that
$$\sum_n\Pr\biggl\{E_Y\Big|\sum_{i,j\le
2^n}\varepsilon_i\pi_2h_n(X_i,Y_j)\Big|\ge C2^n\log n\biggr\}<\infty
\eqno(4.14)$$
for some $C<\infty$,
and that, likewise,
$$\sum_n\Pr\biggl\{E_X\Big|\sum_{i,j\le
2^n}\pi_2h_n(X_i,Y_j)\Big|\ge C2^n\log n\biggr\}<\infty.
\eqno(4.15)$$
Then, (4.13)-(4.15) give (4.5) by Lemma 4.1,
concluding the proof of Step 3.

Steps 1-3 together show that
$$
 \sum_{n}\Pr\biggl\{\Bigl|\sum_{i,j\leq 2^{n}}\pi_2h(X_{i},Y_{j})\Bigr|
 \geq C2^{n}\log n\biggr)<\infty,
$$
concluding the proof of the theorem.
\qed

\vskip.1truein
\n{\bf 5. Arbitrary kernels. Final comments.} We conclude with the proof of
Theorem 1.1, a conjecture on the LIL for kernels of more than two
variables, and
several remarks on the limsup in (1.1) and the limit set of the LIL
sequence.

\vskip.05truein
\n{\smc Proof of Theorem 1.1}. Conditions (1.2) and (1.3) are sufficent for the
LIL for degenerate kernels by Theorems 3.1 and 4.4.

If the kernel $h$ satistifies the LIL (1.1), then it satisfies the
decoupled and
randomized LIL by Lemma 2.1(b). Then, by Theorem 3.1, it also satisfies
conditions (1.2) and (1.3). So, it suffices to prove that if the LIL (1.1)
holds
then the kernel $h$ is canonical.

 Since by (1.2)  $E|\pi_2h|^{p}<\infty$ for any $p<2$, we have by the
 Marcinkiewicz type strong law of large numbers for $U$-statistics (Gin\'e and
Zinn, 1992, theorem 2),
$$\lim_{n\rightarrow \infty}{1\over n^{2/p}}\sum_{i\neq j\leq n}
   \pi_2h(X_{i},Y_{j})=0\ \ {\rm a.s.\ \ for\ all}\ 0<p<2.\eqno(5.1)$$
 The LIL for $h$ implies the decoupled LIL (2.8) by Lemma 2.2(a), and therefore
also that
$$\lim_{n\rightarrow \infty}{1\over n^{2/p}}\sum_{i\neq j\leq n}
   h(X_{i},Y_{j})=0\ \ {\rm a.s.\ \ for\ all}\ 0<p<2.\eqno(5.2)$$
 Subtracting (5.1) from (5.2) and using the Hoeffding decomposition we
obtain
$$\lim_{n\rightarrow \infty}n^{1-2/p}\Bigl|\sum_{i\leq n}\bigl(
   \pi_1h(X_{i})-{1\over 2}Eh\bigr)+\sum_{j\leq n}\bigl(
   \pi_1h(Y_{j})-{1\over2}Eh\bigr)\Bigr|=0\ \ {\rm a.s.}$$
 However if $p\geq 4/3$ this yields, by the  CLT or the
LIL in $\bf R$, that
 $$\pi_1h(X)-{1\over2}Eh=0\ \ {\rm a.s.}$$
Since $\pi_1h$ is centered, it follows that $E h=0$ and $\pi_1h(X)=0$ a.s.
Hence
$h=\pi_2h$ is canonical for the law of $X$. \qed

\vskip.05truein
The following conjecture for kernels of more than two variables seems only
natural.

\proclaim{\smc Conjecture 5.1}.
 Let $h$ be a  kernel of $d$ variables symmetric in its entries. Then $h$
satisfies the law of the iterated
 logarithm
$$\limsup_{n\rightarrow\infty}{1\over(n\log\log n)^{d/2}}\Bigl|
   \sum_{1\leq i_{1}<i_{2}<\ldots<i_{d}\leq n}h(X_{i_{1}},\ldots,X_{i_{d}})
   \Bigr|<\infty \ \  a.s.\eqno(5.3)$$
if and only if the following conditions hold:\hfil\break
a) h  is canonical for the law of $X$, that is
$E_{X_{i}}h(X_{1},\ldots,X_{d})=0$
a.s.\hfil\break and there exists $C<\infty$ such that\hfil\break
b) $$E\min(h^{2},u)\leq C(L_{2}u)^{d-1}\eqno(5.4)$$ for all
$u>0$,
 and\hfil\break
c) $$\sup\bigl\{E[h(X_{1},\ldots,X_{d})\prod_{i=1}^df_{i}(X_{i})]:
   Ef_{i}^{2}(X)\leq 1,\|f_i\|_\infty<\infty,
i=1,\ldots,d\bigr\}<\infty.\eqno(5.5)$$

We know at present that the necessity part of this conjecture is true.

The problem of determining the lim sup in (1.1) when $E h^2=\infty$ is open
and, a fortiori, so is the problem of determining the limit set of the LIL
sequence. We now briefly comment on these questions. The previous results
do give
the order of the limsup in (1.1) up to constants as we show next.
In the theorem that follows we denote the quantity in (1.3) as
$\|h\|_{L_2\mapsto L_2}$.

\proclaim{\smc Theorem 5.2}. Suppose that $h(x,y)$ is canonical for the law of
$X$. Then there exists a universal constant $C$ such that, almost surely,
$$\eqalignno{C^{-1}\biggl[\|h\|_{L_2\mapsto
L_2}&+\limsup_{u\to\infty}\sqrt{{E(h^2\wedge u)\over L_2u}}\biggr]&\cr
&\le
\limsup_{n\to\infty}{1\over nL_2n}\Big|\sum_{1\le i\le j\le
n}h(X_i,X_j)\Big|&\cr &~~~~~~~~~\le C\biggl[\|h\|_{L_2\mapsto
L_2}+\limsup_{u\to\infty}\sqrt{{E(h^2\wedge u)\over L_2u}}\biggr].&(5.6)\cr}$$
The same inequality holds true if $h$ is arbitrary and $h(X_i,X_j)$ is replaced
in (5.6) by the randomized $\ve_i\ve_jh(X_i,X_j)$, or by the decoupled
versions.

\n{\smc Proof}. Lemma 2.1 and the proof of necessity of Theorem 3.1 (see also
Corollary 2.4) give the left hand side bound for decoupled and randomized
kernels. The right hand side bound, also for decoupled and randomized kernels,
follows from the proof of sufficiency of Theorem 3.1: let
$$K:=\max\biggl[\|h\|_{L_2\mapsto
L_2},\limsup_{u\to\infty}\sqrt{{E(h^2\wedge u)\over L_2u}}\biggr];$$
if $K=1$, the proof of Theorem 3.1 produces (3.4) for a fixed constant $C$ that
could be computed if necessary, as can be seen from steps 7
and 8 (the only ones that contribute to the limsup), and if $K\ne 1$, (3.4)
with
$C$ replaced by $CK$
is obtained by considering the kernel
$h/K$.  Then, Corollary 2.5 yields the right hand side of (5.6).
De-randomization as in Section 4 gives the bounds (5.6) for
canonical kernels.
\qed

\vskip.05truein

We
know that when
$Eh^2<\infty$ and
$h$ is a canonical kernel of $d$ variables, the limsup in (5.3) is just the
quantity in (5.5), and even more, that  the limit set of the sequence
$$\biggl\{{d!\over(2n\log\log n)^{d/2}}
   \sum_{1\leq i_{1}<i_{2}<\ldots<i_{d}\leq n}h(X_{i_{1}},\ldots,X_{i_{d}})
   \biggr\}$$
is a.s.
$$\bigl\{E[h(X_{1},\ldots,X_{d})\prod_{i=1}^df(X_{i})]:
   Ef^{2}(X)\leq 1\bigr\}$$
(Dehling,
1989, for $d=2$ and Arcones and Gin\'e, 1995, in general). Then, restricting to
kernels of two variables, several concrete questions arise: 1) is any of
the two
summands in the bounds (5.6) superfluous?; 2) at least in the case when the
kernel $h$ defines a compact operator of $L_2$, can we determine the limit set
of the LIL sequence from the limit set for finite rank $h$ by operator
approximation?, and of course, 3) what is the limit set in general? We will
answer 1) by means of examples showing that, in general, both summands in the
bound (5.6) are essential, and, regarding question 2) we will also
determine the limit set for a class of kernels that induce compact operators in
$L_2$. We wil show, moreover, that there are kernels
$h$ that give non-compact operators for which the LIL holds (the
examples in Gin\'e and Zhang (1996) define compact operators and
suitable modifications will give non-compact ones). Finally, question 3) will
remain open but we will show that the limit set is always an interval.

\vskip.05truein
\n{\smc Example 5.3}. We consider the kernel
$$h(x,y)=\sum_{n=1}^\infty{a_n\over b_n}I_n(x)I_n(y),\eqno(5.7)$$
where $\{I_n\}$ is a sequence of functions on $\bf R$ with disjoint supports
contained in $[0,1]$ such that $\int_{\bf R}I_n(u)du=0$, $I_n(x)\in\{-1,0,1\}$
for each $x\in\bf R$, the sequence $\{b_n\}$ is defined by $b_n=\int_{\bf
R}I_n^2(u)du$ and $\{a_n\}$ is an arbitrary bounded sequence of real numbers.
Then, if, as will be the case, for $X,Y$ i.i.d. uniform on $[0,1]$,
$E|h(X,Y)|<\infty$,
$h$  is a canonical kernel for the uniform distribution on
$[0,1]$. Since $\{b_n^{-1/2}I_n\}$ is an orthonormal sequence in
$L_2:=L_2({\cal L}(X))$, we have
$$\|h\|_{L_2\mapsto L_2}=\sup_{n\in\bf N}|a_n|.\eqno(5.8)$$
If we further assume that $\{a_n/b_n\}$ is an incresing sequence, then
$$\limsup_{u\to\infty}{E(h^2\wedge u)\over
L_2u}=\limsup_n{\sum_{k=1}^na_k^2+{a_n^2\over b_n^2}\bigl(\sum_{k=n+1}^\infty
b_k^2\bigr)\over L_2(b_n^{-1})}.$$ So, if we choose $a_n=a$ for all $n$ and
$I_n$ such that
$b_n=\exp\bigl[-\exp(a^2n/b)\bigr]$ for large $n$, then
$$\limsup_{u\to\infty}{E(h^2\wedge u)\over L_2u}=b.\eqno(5.9)$$
Thus, in this case, the kernel $h$ satisfies the LIL by Theorem 3.1.
Moreover, (5.8) and (5.9) show that the two quantities appearing in the bounds
(5.6) are not comparable (and, in particular, neither of them is superfluous).
In this type of examples, the operator in $L_2$ with kernel $h$ is compact
if and
only if
$\lim_n a_n=0$, thus showing that there are canonical kernels $h$ which satisfy
the LIL but that do not define a compact operator on $L_2$.

If $Eh^2<\infty$, then the operator norm dominates the bound in (5.6), as the
limsup of the normalized truncated second moments of $h$ is zero. Even for
kernels
$h$ defining  compact operators we may have that it is this second term that
dominates the bound: for $a_n=1/\sqrt{n}$ and $b_n=2^{-n}$, consider
the kernels $h_m(x,y)=\sum_{n=m}^\infty a_nb_n^{-1}I_n(x)I_n(y)$; then we have
$\|h_m\|_{L_2\mapsto L_2}=1/\sqrt{m}\to0$ whereas
\hfil\break
$\limsup_{u\to\infty}{E(h_m^2\wedge u)\over L_2u}=1$ for all $m$.

\vskip.05truein
There is, however, a class of canonical kernels $h$ satisfying the LIL
and defining compact operators for which the limit set of the LIL sequence is
the numerical range of the operator defined by $h$, as is the case when $h$ has
finite second moment. In the next proposition  $H$ will denote the operator on
$L_2$ defined by extension  of the equation
$Hf(y)=Eh(X,y)f(X)$, $f\in L_\infty({\cal L}(X))$ (this operator
exists under condition (1.3)).

\proclaim {\smc Proposition 5.4}. Let $h$ be a canonical kernel for
the law of
$X$ such that \hfil\break
a) $$\limsup_{u\to\infty}{E(h^2\wedge u)\over L_2u}=0\eqno(5.10)$$
and\hfil\break
b) the operator $H$ is a compact operator on $L_2({\cal L}(X))$.\hfil\break
Then, the limit set of the sequence
$$\biggl\{{1\over 2nL_2n}\sum_{1\le i\ne j\le n}h(X_i,X_j)\biggr\}\eqno(5.11)$$
is almost surely the closure of the set
$$\Bigl\{E h(X,Y)f(X)f(Y):E f^2(X)\le
1,\|f\|_\infty<\infty\Bigr\},\eqno(5.12)$$
that is, the numerical range of the operator $H$,
$\{E(f(X)Hf(X)):Ef^2(X)\le1\}$.

\n{\smc Proof}. We set, from now on, $L_2:=L_2({\cal L}(X))$. The proof
consists in approximating the operator $H$ with kernel $h$ by suitable
operators $H_m$ with simple kernels, in particular, square integable
kernels. We
begin by showing that there exists an increasing sequence
${\cal G}_m$ of finite sub-$\sigma$-algebras of $\cal S$ such that, if $P_m$
denotes the orthonormal projection onto the subspace of ${\cal G}_m$-measurable
functions,
$$\|P_mHf- Hf\|_{L_2}\to0, \ \ f\in L_2.$$
Indeed, $H$ being a compact operator, its range is a separable set in $L_2$.
Therefore we can find a sequence $\{g_i\}\subset L_2$ of simple functions such
that the range of $H$ is contained in the closure of the sequence $\{g_i\}$.
Now, it is enough to set
$${\cal G}_m:=\sigma(g_1,\dots,g_m)$$
to get the desired property. This is so because, obviously, $P_mg_i\to g_i$ for
each $i\in \bf N$, and the set $\{f\in L_2:P_mf\to f\ {\rm in}\ L_2 \ {\rm
norm}\}$ is closed in view of $\|P_m\|_{L_2\mapsto L_2}\le 1$.

For each $m\in\bf N$ we define
$$h_m(x,y)=\sum_{A,B\ {\srm atoms\ of}\ {\cal G}_m \atop
\Pr\{X\in A,Y \in B\} \ne 0}{Eh(X,Y)I_A(X)I_B(Y)\over \Pr\{X\in A\}\Pr\{Y\in
B\}}I_A(x)I_B(y),$$
where, as usual, $Y$ is an independent copy of $X$. In other words, $h$ is
defined by the condition
$$h_m(X,Y)=E\bigl(h(X,Y)\big|\sigma(X^{-1}({\cal G}_m),Y^{-1}({\cal
G}_m))\bigr).$$
The operator $H_m$ of $L_2$ with kernel $h_m$
satisfies
$H_m=P_mHP_m$, as is seen from its definition.  Then, since
$\|P_mHf-Hf\|_{L_2}\to 0$ for any
$f\in L_2$, and since $H$ is a compact operator in $L_2$, we obtain that
$$\lim_{n\to\infty}\|H_m-H\|_{L_2\mapsto L_2}=0.\eqno(5.13)$$
To see this, we note that, since $(P_m-I)H$ is the adjoint of $H(P_m-I)$ and
$P_m$ has norm 1,
$$\|H_m-H\|_{L_2\mapsto L_2}=\|P_mH(P_m-I)+(P_m-I)H\|_{L_2\mapsto
L_2}\le2\|(P_m-I)H\|_{L_2\mapsto L_2};$$
now
(5.13) follows by a simple compactness argument.

The result follows from the previous observation together with Theorem 5.2
applied to $h_m$ and to $h-h_m$, by a standard approximation argument that we
now sketch. Before we do this, we should note that the closure in $L_2$ of the
set (5.12) is the numerical range of $H$ because bounded functions are dense in
$L_2$, the unit ball of
$L_2$ is weakly compact and if $f_n\to f$ weakly, with
$\|f_n\|_{L_2}\le 1$, then, by compactness of $H$, $Hf_n\to Hf$ weakly. Let us
write $\langle\cdot,\cdot\rangle$ for the inner product in $L_2$, set
$$L:=\{\langle Hf,f\rangle:\|f\|_{L_2}\le 1\}$$
and, for any kernel $g(x,y)$ of two variables,
$$\alpha_n(g):={1\over 2nL_2 n}\sum_{1\le i\ne j\le n}g(X_i,X_j).$$
If $x\in L$ let $f\in L_2$ with $\|f\|_{L_2}\le 1$ be such that
$x=\langle Hf,f\rangle$. Then, by the LIL for kernels with finite second
moment,
given $m\in\bf N$, for almost every $\omega$ there is a subsequence
$n_{k(\omega)}$ such that
$$\alpha_{n_{k(\omega)}}(h_{m}(\omega))\to\langle H_mf,f\rangle.\eqno(5.14)$$
Also, since $h$ satisfies (5.10) and $h_m$ has finite second moment, Theorem
5.2 gives
$$\limsup_n|\alpha_n(h_m-h)|\le K\|H_m-H\|_{L_2\mapsto L_2}\ \ {\rm
a.s.}.\eqno(5.15)$$
Moreover, by (5.13),
$$\langle H_mg,g\rangle\to\langle Hg,g\rangle,\ \ g\in L_2.\eqno(5.16)$$
Combining these three limits we obtain that
$x$ is a.s. a limit point of the sequence $\{\alpha_n(h)\}$. Conversely,
suppose
now that $x$ is a limit point of this sequence. Then, by (5.15), given
$\varepsilon>0$, for all  $m$ large enough and  for almost every $\omega$
there
exists a subsequence
$n_{k(\omega)}$ such that
$$|x-\alpha_{n_{k(\omega)}}(h_{m}(\omega))|<{\varepsilon\over 2}.$$
Therefore, by the LIL for square integrable kernels and (5.16), there is $f\in
L_2$ with $\|f\|_{L_2}\le 1$ such that
$$|x-\langle Hf,f\rangle|<\ve.$$ So, taking $\ve=1/n$, there is a sequence
$f_n$
in the unit ball of $L_2$ such that
$$x=\lim_n\langle Hf_n,f_n\rangle.$$
Since the unit ball of $L_2$ is weakly compact, the sequence $\{f_n\}$ has a
subsequence  $\{f_{n_k}\}$ that converges weakly to a function $f$ in the unit
ball of $L_2$. It then follows by compactness of $H$ that $x=\langle
Hf,f\rangle$, that is, $x\in L$. \qed

\vskip.05truein
For example the previous proposition applies to the kernels $h$ of Example 5.2
for
$a_n=n^{-1/2}\ell(n)$ and $b_n=2^{-n}$, where $\ell(n)$ is any slowly varying
function tending to zero as $n\to\infty$. However, if $\ell(n)=1$ then $h$
still satisfies the LIL (1.1) by Theorem 1.1 and defines a compact operator in
$L_2$, but Proposition 5.4 does not apply to it; actually, we do not know what
the limit set is in this case.

As mentioned, the problem of determining the a.s. limit set of the sequence
(5.11)  in the general case remains open but we can show that it is an
interval.

\proclaim{\smc  Proposition 5.5}. Let $h$ be a canonical kernel satisfying
conditions (1.2) and (1.3). Then, the limit set of the LIL sequence (5.11) is
an interval.

\n{\smc Proof}. To prove that the limit set of the sequence (5.11) is an
interval, it suffices to show that the difference of two consecutive
terms of the sequence
tends to zero a.s.  By
(1.2) and the law of large numbers for $U$-statistics (or by the LIL), this
reduces to showing that
$${1\over n\log\log n}\sum_{1\le i< n}
h(X_{i},X_{n})\to0\ \ {\rm a.s.}\eqno(5.17)$$
We will first prove
$${1\over n\log\log n}\sum_{1\le i< n}
\ve_ih(X_{i},Y_{n})\to0\ \ {\rm a.s.}\eqno(5.18)$$
and then will show that $\ve_i$ can be removed and that $Y_n$ can be replaced
by $X_n$.

To prove (5.18), it is enough to prove that
for all $\delta>0$
$$\sum_{n}\Pr\biggl\{\max_{2^{n-1}< k\le 2^n}{1\over2^n\log{n}}
   \Bigl|\sum_{1\le i< k}\ve_{i}h(X_{i},Y_{k})\Bigr|>\delta
\bigg\}<\infty\eqno(5.19)$$
(see e.g. the proof of Corollary 2.4).
Let $h_{n}=hI_{A_{n}}$ and $\tilde{h}_{n}=h-h_{n}$,
where
$$A_{n}=\bigl\{(x,y):|h(x,y)|\le 2^{n}\log n, f_{n}(y)\le
     2^{n}(\log n)^2\bigr\}.$$
Then as in Steps 1 and 2 of the proof of Theorem 3.1 we get
$$\sum_{n}\Pr\biggl\{\max_{2^{n-1}< k\le 2^n}{1\over2^n\log{n}}
   \Bigl|\sum_{1\le i< k}\ve_{i}\tilde{h}_{n}(X_{i},Y_{k})\Bigr|>\delta \bigg\}
   <\infty.$$
In order to prove
$$\eqalign{\sum_{n}\Pr\biggl\{&\max_{2^{n-1}<k\le
2^n}{1\over2^n\log{n}}\Bigl|
   \sum_{1\le i< k}\ve_{i}h_{n}(X_{i},Y_{k})\Bigr|>\delta \bigg\}\cr
&\le \sum_{n}2^n\Pr\biggl\{\Bigl|\sum_{1\le i< 2^n}
   \ve_{i}h_n(X_{i},Y)\Bigr|>\delta 2^n\log{n}\bigg\}<\infty.}$$
we apply Chebyshev's inequality as in Step 3, reducing the above inequality
to convergence of the two series
$$\sum_n {1\over 2^{2n}(\log{n})^4}Eh_n^4(X,Y)<\infty,$$
$$\sum_n {1\over 2^n(\log{n})^4}Eh_n^2(X_1,Y)h_n^2(X_2,Y)<\infty.$$
But these two series converge, just like the first
and second series in Step 3. (5.19) is thus proved.

Next we show that we can remove the Rademacher variables from (5.18), that is,
that (5.18) implies
$${1\over n\log\log n}\sum_{1\le i< n}
h(X_{i},Y_{n})\to0\ \ {\rm a.s.}\eqno(5.20)$$
Let $\{\tilde{X_i}\}$ be a copy of $\{X_i\}$, independent of
$\{X_i\}$ and $\{Y_i\}$, and set
$$\xi_n:={1\over n\log\log n}\sum_{1\le i< n}
h(X_{i},Y_{n}),\ \ \tilde{\xi}_n:={1\over n\log\log n}\sum_{1\le i< n}
h(\tilde{X}_{i},Y_{n}).$$
If (5.18) holds, then $\xi_n-\tilde{\xi}_n\to0$
a.s. by Fubini's theorem and the equidistribution of the variables $X_i$.
Hence, (5.20) will follow by a standard argument if $\xi_n\to0$ in probability
conditionally on the sequence $\{Y_i\}$. So, assuming (wlog) that the variables
$X$ and $Y$ are defined on different factors of a product probability space
$\Omega'\times\Omega$, we must show that
$${1\over a_n}\sum_{1\le i<n}h(X_i,Y_n(\omega))\to0\ \ {\rm in\ pr.},\
\omega-{\rm a.s.},\eqno(5.21)$$
where, for ease of notation, we set $a_n:=(nL_2 n)^{-1}$. Now, since
$${1\over a_n}
\sum_{1\le i< n}
\ve_ih(X_{i},Y_{n})\to 0\ \ {\rm in\ pr.},\ \omega-{\rm a.s.}$$
by (5.18),
L\'evy's inequality applied conditionally on $\{Y_i\}$ gives
$$n{\Pr}_X\bigl\{|h(X,Y_n)|>a_n\bigr\}\to 0 \hbox{ a.s. }\eqno(5.22)$$
and then, Hoffmann-J\o rgensen's inequality  applied conditionally
after truncation, yields
$${n\over a_n^2}E_Xh^2(X,Y_n)I_{\{|h(X,Y_n)|\le a_n\}}\to 0 \hbox{
a.s.}\eqno(5.23)$$
Moreover,
$${n\over a_n}E_Xh(X,Y_n)I_{\{|h(X,Y_n)|\le a_n\}}\to 0 \hbox{
a.s.}\eqno(5.24)$$
To prove that this last limit holds, note first that, since $E_Xh=0$,
$$E_Xh(X,Y_n)I_{|h(X,Y_n)|\le a_n}=E_Xh(X,Y_n)I_{\{|h(X,Y_n)|> a_n\}},$$
and then that
$$\sum_n{n\over a_n}E|h(X,Y)|I_{|h(X,Y)| >a_n} <\infty$$
because, after exchanging expectation and sum and then summing on $n$, we
see that this series is bounded by a constant times $E{h^2\over L_2^2|h|}$,
which is finite. Now, (5.22)-(5.24) give that, for all $\ve>0$,
$$\eqalign{{\Pr}_X\biggl\{{1\over a_n}\Big|\sum_{1\le i\le n}
h(X_{i},Y_{n})\Big|>\ve\biggr\}&\le  n{\Pr}_X\bigl\{|h|>a_n\bigr\}+
I_{\{na_n^{-1}|E_XhI_{\{|h|\le a_n\}}|>\ve/2\}}\cr
&~~~~~~~~~~~~~~~~~~~~~~~+{4\over\ve^2}{n\over
a_n^2}E_Xh^2I_{\{|h|\le a_n\}}\to0\ \ {\rm a.s.},\cr}$$
proving (5.21), hence, (5.20).

Finally, to undecouple, assume (5.20) holds.
By Theorem 1.1 and the $0-1$ law we know that
$$\limsup_n{1\over n\log\log n}\Big|\sum_{1\le i< n}h(X_i,X_n)\Big|=C\ \ {\rm
a.s.}\eqno(5.25)$$
for some $C<\infty$, and must show that $C=0$.
Then, we can assume that this limsup is attained by the sequence of even terms,
that is,
$$\limsup_n{\big|\sum_{1\le i< 2n}h(X_i,X_{2n})\big|\over 2n\log\log (2n)}=C\ \
{\rm a.s.}\eqno(5.26)$$
(otherwise we can take the subsequence of odd terms from (5.25) and continue in
the same way as we will now proceed). But
$$\eqalign{\limsup_n{1\over
2n\log\log (2n)}&\Big|\sum_{1\le i< 2n}h(X_i,X_{2n})\Big|\cr
&\le
\limsup_n{1\over 2n\log\log (2n)}\Big|\sum_{1< i< 2n\atop i\ {\srm
even}}h(X_i,X_{2n})\Big|\cr
&~~~~~~~~~~+\limsup_n{1\over
2n\log\log (2n)}\Big|\sum_{1\le i< 2n\atop i\ {\srm
odd}}h(X_i,X_{2n})\Big|\cr
&=\limsup_n{1\over2n\log\log(2n)}\Big|\sum_{1\le i<n}h(X_i,X_n)\Big|\cr
&~~~~~~~~~~+\limsup_n{1\over2n\log\log(2n)}\Big|\sum_{1\le
i<n+1}h(X_i,Y_{n+1})\Big|\cr
&= {C\over 2}\cr}$$ by (5.25) and (5.20). This
contradicts (5.26) unless $C=0$, proving (5.17).
\qed

\vskip.3truein

\n{\bf Acknowledgements.} The second and third named authors
carried out the research for the present article respectively at the
departments
of Mathematics of Texas A\&M University and Georgia Tech, and wish to
acknowledge their hospitality.

\vskip.4truein
\centerline{\bf References}
\vskip.2truein \baselineskip=10pt
\parskip=4pt

{\smc\n Arcones, M. and Gin\'e, E.} (1995). On the law of the
iterated logarithm for canonical $U$--statistics and processes.
{\it Stoch. Proc. Appl.} {\bf 58} 217-245.\par
{\smc\n Dehling, H.} (1989). Complete convergence of
triangular arrays and the law of the iterated logarithm for
degenerate $U$--statistics. {\it Stat. Probab. Letters} {\bf 7}
319--321.\par
{\smc\n Dehling, H.; Denker, M. and Philipp, W.} (1984). Invariance
principles for von Mises and $U$-statistics. {\it Zeits. Wahrsch. verw. Geb.}
{\bf 67} 139-167.\par
{\smc\n Dehling, H.; Denker, M. and Philipp, W.}  (1986) A bounded
law of the iterated logarithm for Hilbert space valued martingales
and its application to $U$--statistics. {\it Prob. Th. Rel.
Fields} {\bf 72} 111-131.\par
{\smc\n de la Pe\~na, V. and Montgomery--Smith, S.} (1994).
Bounds for the tail probabilities of $U$-statistics and quadratic forms.
{\it Bull. Amer. Math. Soc.} {\bf 31} 223-227.\par
{\smc\n Gin\'e, E. and Zhang, C.-H.} (1996). On integrability in
the LIL for degenerate $U$--statistics. {\it J. Theoret. Probab.} {\bf
 9} 385--412.\par
{\smc\n Gin\'e, E. and Zinn, J.} (1992). Marcinkiewicz type laws
of large numbers and convergence of moments for $U$--statistics.
{\it Probability in Banach Spaces 8} 273--291. Birkh\"auser,
Boston.\par
{\smc\n Gin\'e, E. and Zinn, J.} (1994). A remark on convergence
in distribution of $U$--statistics. {\it Ann. Probab.} {\bf
22} 117--125.\par
{\smc\n Goodman, V.} (1996). A bounded LIL for second order
$U$--statistics. Preprint.\par
{\smc\n Halmos, P. R.} (1946). The theory of unbiased
estimation. {\it Ann. Math. Statist.} {\bf 17} 34--43.\par
%{\smc\n Hitczenko, P.} (1988). Comparison of moments for tangent sequences
%of random variables. {\it Probab. Th. Rel. Fields} {\bf 78} 223-230.\par
{\smc\n Hoeffding, W.} (1948). A class of statistics with
asymptotically
normal distribution. {\it  Ann. Math. Statist.}  {\bf 19}
293-325.\par
{\smc\n Lata\l a, R.} (1999). Tails and moment estimates for some type
of chaos. {\it Studia Math.}, to appear.\par
{\smc\n Lata\l a, R. and Zinn, J.} (1999).
Necessary and sufficient conditions for the strong law of large numbers for
$U$-statistics. Preprint.\par
{\smc\n Ledoux, M.} (1996). On Talagrand's deviation inequalities for product
measures. {ESAIM, P\&S},{\bf 1} 63-87. (http://www.emath.fr/Maths/Ps)\par
{\smc \n Montgomery--Smith, S.} (1993). Comparison of sums of
independent identically dsitributed random variables. {\it
Prob. Math. Statist.}  {\bf 14} 281--285.\par
{\smc\n Rubin, M. and Vitale, R. A. } (1980).
Asymptotic distribution of symmetric statistics. {\it Ann.  Statist.} {\bf 8}
165-170.\par
{\smc\n Serfling, R. J.}  (1971). The law of the iterated logarithm
for $U$--statistics and related von Mises functionals. {\it Ann.
Math. Statist.} {\bf 42} 1794.
\par
{\smc\n Talagrand, M.} (1994). Sharper bounds for Gaussian and empirical
processes. {\it Ann. Probab.} {\bf 22} 28-76.

{\smc\n Talagrand, M.} (1996). New concentration inequalities in product
spaces. {\it Invent. Math.} {\bf 126} 505-563.\par
\n{\smc Teicher, H.} (1995). Moments of randomly stopped sums revisited.
{\it J.
Theoret. Probab.} {\bf 8} 779-794.\par
{\smc \n Zhang, C.-H.} (1999). Sub-Bernoulli functions, moment
inequalities and strong laws for nonnegative and symmetrized $U$-statistics.
{\it Ann. Probab.} {\bf 27} 432-453.\par

\vskip.5truein\n
{\srm
\settabs \+Miguel A. Arcones$\phantom{\displaystyle
{la mare que et va
parir tralalalalalalalallallla
ima}}$& Joel Zinn\cr

\+Department of Mathematics&Institute of Mathematics\cr
\+and Department of Statistics&Warsaw University\cr
\+University of Connecticut&Banacha 2\cr
\+Storrs, CT 06269&02-097 Warszawa\cr
\+USA&Poland\cr
\+gine@uconnvm.uconn.edu&kwapstan@mimuw.edu.pl\cr
\+&rlatala@mimuw.edu.pl\cr}

\vskip.2truein

{\srm
\settabs \+aaaaaaaaaaaaaaaaaaaaaaaaaaaaaaaaaaaa&Department of
mathematical Sciences\cr
\+&Department of Mathematics\cr
\+&Texas A\&M University\cr
\+&College Station, TX 77843\cr
\+&jzinn@math.tamu.edu\cr}

\end